\documentclass[number,citesort,seceqn,dvips]{arxbj}
\usepackage{mathbh,dcolumn}
\usepackage{graphicx}

\aid{0}
\volume{16}
\issue{4}
\pubyear{2010}
\firstpage{1039}
\lastpage{1063}
\doi{10.3150/10-BEJ256}

\makeatletter
\newcolumntype{d}[1]{D{.}{.}{#1}}
\def\eqref#1{(\ref{#1})}
\newtheorem{Lemma}{Lemma}
\newremark{Remark}{Remark}[section]
\newtheorem{Theorem}{Theorem}[section]
\newtheorem{Corollary}{Corollary}[section]
\newremark{Example}{Example}[section]

\makeatother

\begin{document}
\begin{frontmatter}

\title{Frontier estimation and extreme value theory}
\runtitle{Monotone frontier estimation}

\begin{aug}
\author[1]{\fnms{Abdelaati} \snm{Daouia}\corref{}\thanksref{1,e1}\ead[label=e1,mark]{daouia@cict.fr}},
\author[1]{\fnms{Jean-Pierre} \snm{Florens}\thanksref{1,e2}\ead[label=e2,mark]{florens@cict.fr}} \and
\author[2]{\fnms{L\'eopold} \snm{Simar}\thanksref{2}\ead[label=e3]{leopold.simar@uclouvain.be}}
\runauthor{A. Daouia, J.-P. Florens and L. Simar}
\address[1]{Toulouse School of Economics (GREMAQ, University of
Toulouse), Manufacture des
Tabacs, 21 All\'{e}e de Brienne, 31000 Toulouse, France.
E-mails: \printead*{e1}, \printead*{e2}}
\address[2]{Institute of Statistics, Catholic University of Louvain,
Voie du Roman Pays 20, B-1348, Louvain-la-Neuve, Belgium.
\printead{e3}}
\pdfauthor{Abdelaati Daouia, Jean-Pierre Florens, Leopold Simar}
\end{aug}

\received{\smonth{2} \syear{2008}}
\revised{\smonth{1} \syear{2010}}

%
\begin{abstract}
In this paper, we investigate the problem of nonparametric monotone frontier
estimation from the perspective of extreme value theory. This enables
us to revisit
the asymptotic theory of the popular free disposal hull estimator in a
more general
setting, to derive new and asymptotically Gaussian estimators and to provide
useful asymptotic confidence bands for the monotone boundary function. The
finite-sample behavior of the suggested estimators is explored via
Monte Carlo experiments. We also apply our approach to a real data set
based on
the production activity of the French postal services.
\end{abstract}

%

\end{frontmatter}
%

\section{Introduction}\label{sec:intr}

In production theory and efficiency analysis, there is sometimes the
need to estimate the boundary of a production set (the
set of feasible combinations of inputs and outputs). This boundary (the
production frontier) represents the set of optimal production plans so that
the efficiency of a production unit (a firm, for example) is obtained
by measuring
the distance from this unit to the estimated production frontier. Parametric
approaches rely on parametric models for the frontier and the underlying
stochastic process, whereas nonparametric approaches offer much more flexible
models for the data-generating process (see, for example, \cite
{DSBOOK} for recent surveys on this topic).

Formally, in this paper, we consider technologies where $x \in\mathbb{R}
^p_+$, a vector of production factors (inputs) is used to produce a
single quantity (output) $y\in\mathbb{R}_+$. The attainable
production set is
then defined, in standard microeconomic theory, as
$ \mathbb{T}= \{(x,y) \in\mathbb{R}^p_+ \times\mathbb{R}_+ \mid x
\mbox{ can produce } y\}.$
Assumptions are usually made on this set, such as free disposability of inputs
and outputs, meaning that if $(x,y)\in\mathbb{T}$, then $(x',y')\in
\mathbb{T}$ for any
$(x',y')$ such that $x'\ge x$ (this inequality must be understood
componentwise) and $y'\le y$.
To the extent that the efficiency of a firm is a concern, the boundary
of~$\mathbb{T}$ is of interest. The efficient boundary (or \textit{production
frontier}) of~$\mathbb{T}$ is the locus of optimal production plans (maximal
achievable output for a given level of inputs). In our setup, the
production frontier is represented by the graph of the production function
$\phi(x) = \sup\{ y \mid  (x,y)\in\mathbb{T}\}.$
The economic efficiency score of a firm operating at the level $(x,y)$
is then given by the ratio $\phi(x)/y$.

Cazals \textit{et al.}~\cite{CAZ} proposed a probabilistic interpretation
of the production frontier.
Let $\mathbb{T}$ be the support of the joint distribution of a random vector
$(X,Y) \in\mathbb{R}^p_+ \times
\mathbb{R}_+$ and let $(\Omega,\mathcal{A},\mathbb{P})$ be the probability
space on which the vector of inputs $X$ and the output $Y$ are defined.
The distribution function of $(X,Y)$ can be denoted $F(x,y)$ and
$F(\cdot|x)=F(x,\cdot)/F_X(x)$ will be
used to denote the conditional distribution function of $Y$ given
$X\leq x$,
with $F_X(x)=F(x,\infty)>0$. It has been proven in \cite{CAZ} that
\begin{eqnarray*}
\varphi(x) = \sup\{y \geq0\mid F(y|x)<1\}
\end{eqnarray*}
is a monotone non-decreasing function with $x$. So, for all $x'\ge x$
with respect to the partial order, $\varphi(x') \ge\varphi(x)$. The
graph of $\varphi$ is the smallest non-decreasing surface which is
greater than or equal to the upper boundary of $\mathbb{T}$. Further,
it has
been shown that under the free disposability assumption, $\varphi
\equiv\phi$, that is, the graph of $\varphi$ coincides with the
production frontier.

Since $\mathbb{T}$ is unknown, it must be estimated from a sample of
i.i.d.~firms
${\mathcal{X}}_n=\{(X_i,Y_i)\mid i=1,\ldots,n\}$. The \textit{free
disposal hull}
(FDH) $\widehat
{\mathbb{T}}_{\mathrm{FDH}} = \{ (x,y)\in\mathbb{R}_+^{p+1} \mid  y \le
Y_i,  x \ge
X_i,
i=1,\ldots,n \}$ of
${\mathcal{X}}_n$ was introduced by \cite{DST}.
The resulting FDH estimator of $\varphi(x)$ is
\begin{eqnarray*}
\hat{\varphi}_1(x) =\sup\{y\geq0\mid \hat{F}(y|x)<1\}=\max_{i\dvtx  X_i
\leq x} Y_i ,
\end{eqnarray*}
where
$\hat{F}(y|x) = \hat{F}_n(x,y)/\hat{F}_{X}(x)$ with
$\hat{F}_n(x,y)=(1/n) \sum^n_{i=1} \mathbh{1}(X_i\leq x,Y_i\leq y)$ and
$\hat{F}_X(x)=\hat{F}_n(x,\infty)$.
This estimator represents the lowest
monotone step function covering all of the data points $(X_i,Y_i)$.
The asymptotic behavior of $\hat{\varphi}_1(x)$ was first derived
by \cite{KST} for the consistency and by \cite{PaSW,HWA} for the
asymptotic sampling distribution. To summarize, under
regularity conditions, the FDH estimator $\hat{\varphi}_1(x)$ is consistent
and converges to a Weibull distribution with some unknown parameters.
In Park \textit{et al.}~\cite{PaSW}, the
obtained convergence rate $n^{-1/(p+1)}$ requires that the joint
density of $(X,Y)$ has
a jump at its support boundary. In addition, the estimation of the
parameters of the Weibull distribution requires the specification of
smoothing parameters
and the resulting procedure has very poor accuracy. In Hwang \textit{et
al.}~\cite{HWA}, the
convergence of $\hat{\varphi}_1(x)$ to the Weibull distribution was
established in a general case where the
density of $(X,Y)$ may decrease to zero or increase toward infinity at
a speed of
power $\beta$
($\beta>-1$) of the distance from the frontier.
They obtain the convergence rate $n^{-1/(\beta+2)}$ and extend the particular
result of Park \textit{et al.}~\cite{PaSW} where $\beta=0$, but their
result is only derived
in the simple case of one-dimensional inputs $(p=1),$ which may be of less
interest in practice.

In this paper, we first analyze the properties of the FDH estimator
from an
extreme value theory perspective. In doing so, we generalize and extend the
results of Park \textit{et al.}~\cite{PaSW} and Hwang \textit{et
al.}~\cite{HWA} in at least three
directions. First, we provide the necessary and sufficient condition
for the FDH
estimator to converge in distribution and we specify
the asymptotic distribution with the appropriate rate of
convergence. We also provide a limit theorem for moments in a general
framework. Second, we show how the unknown parameter $\rho_x>0,$
involved in
the necessary and sufficient extreme value conditions, is linked to the
dimension $p+1$ of the data and to the shape parameter $\beta>-1$ of
the joint
density: in the general setting
where $p\geq1$
and $\beta=\beta_x$ may depend on $x$, we obtain, under a convenient
regularity
condition, the general convergence rate $n^{-1/\rho_x}=n^{-1/(\beta
_x+p+1)}$ of
the FDH estimator $\hat{\varphi}_1(x)$.
Third, we suggest a strongly consistent and asymptotically normal
estimator of
the unknown parameter $\rho_x$ of the asymptotic Weibull distribution of
$\hat{\varphi}_1(x)$. This also answers the important question of how to
estimate the shape parameter $\beta_x$ of the joint
density of $(X,Y)$ when it approaches the frontier of the support
$\mathbb{T}$.

By construction, the FDH estimator is very non-robust to extremes. Recently,
Aragon \textit{et al.}~\cite{ARA} constructed an original estimator of
$\varphi(x)$, which is more robust than $\hat{\varphi}_1(x),$ but
which keeps the
same limiting Weibull distribution as $\hat{\varphi}_1(x)$ under the
restrictive condition $\beta=0$.
In this paper, we provide further insights and generalize their main
result. We also
suggest attractive estimators of $\varphi(x)$
converging to a normal distribution, which appear to be robust to outliers.
The paper is organized as follows. Section~\ref{sec2} presents
the main results of the paper. Section~\ref{sec3} illustrates how the theoretical
asymptotic results behave in finite-sample situations and gives an example
with a real data set on the production activity of the French postal
services. Section~\ref{sec4} concludes the paper, with proofs deferred for the \hyperref[append]{Appendix}.

\section{The main results}
\label{sec2}

From now on, we assume that $x\in\mathbb{R}^{p}_+$ such that
$F_X(x)>0$ and
will denote by $\varphi_{\alpha} (x)$ and
$\hat{\varphi}_{\alpha}(x)$, respectively, the $\alpha$-quantiles
of the distribution function $F(\cdot|x)$ and its empirical version
$\hat{F}(\cdot|x)$,
\[
\varphi_{\alpha} (x) = \inf\{y \geq0\mid F(y|x) \geq\alpha\}
\quad \mbox{and}\quad  \hat{\varphi}_{\alpha}(x) = \inf
\{y\geq0\mid \hat{F}(y|x) \geq\alpha\}
\]
with $\alpha\in\,]0,1]$. When $\alpha\uparrow1$, the conditional
quantile $\varphi_{\alpha}(x)$ tends to
$\varphi_{1}(x)$, which coincides with the frontier function $\varphi
(x)$. Likewise,
$\hat{\varphi}_{\alpha}(x)$ tends to the FDH estimator
$\hat{\varphi}_{1}(x)$ of $\varphi(x)$ as $\alpha\uparrow1$.

\subsection{Asymptotic Weibull distribution}

We first derive the following interesting results on the problem of
convergence in
distribution of suitably normalized maxima $b^{-1}_n(\hat{\varphi}_{1}
(x) - \varphi(x))$. We will denote by $\Gamma(\cdot)$ the gamma function.

\begin{Theorem}
\label{thm1}
\textup{(i)} If there exist $b_n >0$ and some non-degenerate
distribution function $G_x$ such that
%
\begin{equation}
\label{(2.1)}
b^{-1}_n \bigl(\hat{\varphi}_1 (x) - \varphi(x)\bigr)\stackrel
{{d}}{\longrightarrow}
G_x  ,
\end{equation}
then $G_x(y)$ coincides with $\Psi_{\rho_x}(y)= \exp\{-(-y)^{\rho
_x}\}$ with
support $]{-}\infty, 0]$ for some $\rho_x >0$.
\begin{enumerate}[(iii)]
\item[(ii)] There exists $b_n >0$ such that $b^{-1}_n
(\hat{\varphi}_1 (x) - \varphi(x))$ converges in distribution if and
only if
%
\begin{equation}
\label{(2.2)}
\hspace*{-6pt}\lim_{t\to\infty} \bigl\{1-F\bigl(\varphi(x) -1/tz\mid x\bigr)\bigr\}/\bigl\{1-F\bigl(\varphi(x) -
1/t\mid x\bigr)\bigr\} =
z^{-\rho_x} \qquad  \mbox{for all }  z>0\hspace*{6pt}
\end{equation}
$(\mbox{regular variation with exponent}  -\rho_x,  \mbox{
notation }  1-F(\varphi(x) - \frac{1}{t}\mid x)\in \mathit{RV} _{-\rho
_x})$.

\noindent In this case, the norming constants $b_n$ can be chosen as
$
b_n = \varphi(x) - \varphi_{1-(1/nF_{X}(x))} (x).
$
\item[(iii)] Given \eqref{(2.2)},
$\lim_{n\rightarrow\infty}\mathbb{E}\{b^{-1}_n(\varphi(x)-\hat
{\varphi
}_1 (x)
)\}^k=\Gamma(1+k\rho^{-1}_x)$ for all integers $k\geq1$ and
\begin{eqnarray*}
&&\lim_{n\rightarrow\infty}\mathbb{P}\biggl[\frac{\hat{\varphi}_1
(x) - \mathbb{E}
(\hat{\varphi}_1 (x))
}{\{\operatorname{Var}(\hat{\varphi}_1 (x))\}^{1/2}}\leq y\biggr]\\
&&\quad =\Psi
_{\rho_x}[\{\Gamma(1+2\rho^{-1}_x)-\Gamma^2(1+\rho^{-1}_x)\}
^{1/2}y-\Gamma(1+\rho^{-1}_x)] .
\end{eqnarray*}
\end{enumerate}
\end{Theorem}

\begin{Remark}
\label{rmk2.4}
Since the function $t\mapsto F_X(x)[1 - F(\varphi(x)- \frac{1}{t}
\mid   x)] \in
 \mathit{RV} _{-\rho_x}$ (regularly varying in $t\to\infty$) by
\eqref{(2.2)}, this function can be represented as $t^{-\rho
_x}L_x(t)$ with
$L_x(\cdot) \in \mathit{RV} _{0}$ ($L_x$ being slowly varying) and so
the extreme value condition \eqref{(2.2)} holds if and only if we have
the following representation:
%
\begin{eqnarray}
\label{fxy}
F_X(x)[1 - F(y \vert  x)]=L_x\bigl(\{\varphi(x)-y\}^{-1}\bigr)
\bigl(\varphi(x)-y\bigr)^{\rho_x}\qquad  \mbox{as }  y\uparrow\varphi(x) .
\end{eqnarray}
In the particular case where $L_x(\{\varphi(x)-y\}^{-1}
)=\ell_x$
is a strictly positive function in $x$, it is shown in the next corollary
that $b_n\sim( n\ell_x)^{-1/\rho_x}$.
From now on, a random variable $W$ is said to follow the distribution
$\operatorname{Weibull}(1,\rho_x)$ if $W^{\rho_x}$ is
exponential with parameter $1$.
\end{Remark}

\begin{Corollary}\label{(corcor)}
Given \eqref{fxy} or, equivalently,
\eqref{(2.2)} with $L_x(\{\varphi(x)-y\}^{-1})=\ell
_x>0$, we have
\[
( n\ell_x)^{1/\rho_x}\bigl(\varphi(x) -\hat\varphi_1(x) \bigr)
\stackrel{{d}}{\longrightarrow} \operatorname{Weibull}(1,\rho_x)\qquad  \mbox{as }
n\to
\infty .
\]
\end{Corollary}

\begin{Remark}
\label{rmk2.4.1}
Park \textit{et al.}~\cite{PaSW} and Hwang \textit{et al.}~\cite{HWA}
have obtained similar
results under more restrictive conditions.
Indeed, a unified formulation of the assumptions
used in \cite{PaSW,HWA} can be expressed as
%
\begin{eqnarray}
\label{eq:th4b}
f(x,y) = c_{x}  \{\varphi(x) - y\}^{\beta}+\mathrm{o}\bigl(\{\varphi(x) - y\}
^{\beta}\bigr) \qquad \mbox{as }  y\uparrow\varphi(x),
\end{eqnarray}
where $f(x,y)$ is the joint density of $(X,Y)$, $\beta$ is a
constant satisfying $\beta>-1$ and $c_x$ is a strictly positive
function in
$x$.
Under the restrictive
condition that $f$ is strictly positive on the
frontier (that is, $\beta=0$), Park \textit{et al.}~\cite{PaSW}, among
others, have obtained the limiting Weibull
distribution of the FDH estimator with the convergence rate $n^{-1/(p+1)}$.
When $\beta$ may be non-null, Hwang \textit{et al.}~\cite{HWA} have
obtained the asymptotic Weibull
distribution with the convergence rate $n^{-1/(\beta+2)}$ in the
simple case $p=1$ (here, it is also assumed that \eqref{eq:th4b} holds
uniformly in a
neighborhood of the point at which we want to estimate $\varphi(\cdot
)$, and
that this frontier function is strictly increasing in that
neighborhood and satisfies a Lipschitz condition of order $1$). In the general
setting where $p\geq1$ and $\beta=\beta_x>-1$ may depend on $x$, we have
the following, more general, result, which involves the link between
the tail index $\rho_x$, the data dimension $p+1$ and the shape
parameter $\beta_x$ of the joint density near the boundary.
\end{Remark}

\begin{Corollary}\label{(corscors)}
If the condition of Corollary \ref{(corcor)} holds with $F(x,y)$ being
differentiable near the frontier (that is, $\ell_x>0$,
$\rho_x>p$ and $\varphi(x)$ are differentiable in $x$ with first partial
derivatives of $\varphi(x)$ being strictly positive), then \eqref
{eq:th4b} holds with
$\beta=\beta_x=\rho_x-(p+1)$ and we have
\[
( n\ell_x)^{1/(\beta_x+p+1)}\bigl(\varphi(x) -\hat\varphi_1(x)
\bigr)
\stackrel{{d}}{\longrightarrow} \operatorname{Weibull}(1,\beta_x+p+1)
\qquad \mbox{as
}  n\to
\infty .
\]
\end{Corollary}

\begin{Remark}
\label{rmk2.4.2}
We assume the differentiability of the functions $\ell_x$, $\rho_x$ with
$\rho_x>p$ and $\varphi(x)$ in order to ensure the existence of the joint
density near its support boundary. We distinguish between three different
behaviors of this density at the frontier point
$(x,\varphi(x))\in\mathbb{R}^{p+1}$ based on how the value of $\rho
_x$ compares
to the dimension $(p+1)$: when $\rho_x>p+1$, the joint density decays
to zero
at a speed of power $\rho_x-(p+1)$ of the distance from the frontier; when
$\rho_x=p+1$, the density has a sudden jump at the frontier; when
$\rho_x<p+1$,
the density increases toward infinity at a speed of power $\rho
_x-(p+1)$ of the
distance from the frontier. The case $\rho_x\leq p+1$ corresponds to
sharp or
fault-type frontiers.
\end{Remark}

\begin{Remark}
\label{rmk2.4.22}
As an immediate consequence of Corollary
\ref{(corscors)}, when $p=1$ and $\beta_x=\beta$ (or, equivalently,
$\rho_x=\rho$) does not depend on $x$,
we obtain
the convergence in distribution
of the FDH estimator, as in Hwang \textit{et al.}~\cite{HWA} (see Remark
\ref{rmk2.4.1}), with the same convergence rate $n^{-1/(\beta+2)}$ (in
the notation of \cite{HWA}, Theorem 1,
$\mu(x)=\ell_x(\beta+2)\varphi'(x)=\ell_x\rho_x\varphi'(x)$).
In the other particular case where the joint density is strictly positive
on the frontier, we achieve the best rate of convergence
$n^{-1/(p+1)}$, as in
Park \textit{et al.}~\cite{PaSW} (in the notation of Theorem 3.1 in
\cite{PaSW}, $\mu_{NW,0}/y=\ell_x^{1/(p+1)}=\ell_x^{1/\rho_x}$).

Note, also, that the condition \eqref{eq:th4b} with $\beta=\beta
_x>-1$ (as in Corollary
\ref{(corscors)}) has been considered by \cite{HAR,HAL,Gij}. In
Section \ref{parag2.3}, we answer the important question of how to
estimate the shape
parameter $\beta_x$ in \eqref{eq:th4b} or, equivalently, the regular
variation exponent $\rho_x$ in \eqref{(2.2)}.

As an immediate consequence of Theorem \ref{thm1}(iii) in
conjunction with Corollary \ref{(corscors)}, we obtain
\begin{eqnarray}\label{cod}
\mathbb{E}\{\varphi(x)-\hat{\varphi}_1
(x)\}^k&=&k\{\beta_x+p+1\}^{-1}\{n\ell_x\}^{-k/(\beta_x+p+1)}\Gamma
(k\{\beta_x+p+1\}^{-1})\nonumber\\[-8pt]\\[-8pt]
&&{}+\mathrm{o}\bigl(n^{-k/(\beta_x+p+1)}\bigr) .\nonumber
\end{eqnarray}
This extends the limit theorem of moments of Park \textit{et al.}~(\cite
{PaSW}, Theorem
3.3) to the more general setting where $\beta_x$ may be non-null. Likewise,
Hwang \textit{et al.}~(\cite{HWA}, Remark 1) provide \eqref{cod}
only for $k\in\{1,2\}$, $p=1$ and $\beta_x=\beta$.
The result \eqref{cod} also reflects the well-known curse of
dimensionality from which the
FDH estimator $\hat{\varphi}_1(x)$ suffers as the number $p$ of
inputs-usage increases,
as pointed out earlier by Park \textit{et al.}~\cite{PaSW} in the
particular case where $\beta_x=0$.
\end{Remark}

\subsection{Robust frontier estimators}

By an appropriate choice of $\alpha$ as a function of $n$, Aragon
\textit{et al.}~\cite{ARA} have shown that $\hat{\varphi}_{\alpha
}(x)$ estimates the full frontier $\varphi(x)$
itself and converges to the same Weibull distribution as the FDH
$\hat{\varphi}_1(x)$ under the restrictive conditions of \cite{PaSW}.
The next theorem provides further insights and generalizes their main result.

\begin{Theorem}
\label{Thm2.2}
\begin{enumerate}[(ii)]
\item[(i)] If $b^{-1}_n (\hat{\varphi}_1 (x) -
\varphi(x))\stackrel{{d}}{\longrightarrow} G_x $, then for any fixed integer
$k\geq0$,
\[
b^{-1}_n
\bigl(\hat{\varphi}_{1-k/(n\hat{F}_{X}(x))}(x) - \varphi(x)\bigr)
\stackrel{{d}}{\longrightarrow} H_x \qquad  \mbox{as }  n\to
\infty
\]
for the distribution function $H_x (y)= G_x (y)
\sum^k_{i=0} (-\log G_x (y))^i/i!$.
\item[(ii)] Suppose that the upper bound of the support of $Y$ is finite.
If $b^{-1}_n (\hat{\varphi}_1 (x) - \varphi(x))\stackrel
{{d}}{\longrightarrow } G_x $,
then $b^{-1}_n (\hat{\varphi}_{\alpha_n} (x) -
\varphi(x))\stackrel{{d}}{\longrightarrow} G_x $ for all sequences
$\alpha
_n\rightarrow
1$ satisfying $nb^{-1}_n(1-\alpha_n)\rightarrow0$.
\end{enumerate}
\end{Theorem}

\begin{Remark}
When $\hat{\varphi}_1(x)$ converges in distribution,
the estimator $\hat{\varphi}_{\alpha_n}(x)$,
for $\alpha_n:=1-k/n\hat{F}_X(x)<1$ (that is, $k=1,2,\ldots,$ in Theorem
\ref{Thm2.2}(i)),
estimates $\varphi(x)$ itself and also converges in
distribution, with the same scaling, but a different limit
distribution (here, \mbox{$nb^{-1}_n(1-\alpha_n)\stackrel{\mathrm
{a.s.}}{\longrightarrow}\infty
$)}. To recover the
same limit distribution as the FDH estimator, it suffices to require
that $\alpha_n\to
1$ rapidly so that $nb^{-1}_n(1-\alpha_n)\rightarrow0$.
This extends the main result of Aragon \textit{et al.}~(\cite{ARA},
Theorem~4.3), where the
convergence rate achieves $n^{-1/(p+1)}$ under the restrictive
assumption that the density of $(X,Y)$ is
strictly positive on the frontier.
Note, also, that the estimate $\hat{\varphi}_{\alpha_n}$ does
not envelop all of the data points providing a robust alternative to
the FDH frontier $\hat{\varphi}_1$; see \cite{DNN} for an
analysis of its quantitative and qualitative robustness properties.
\end{Remark}

\subsection{Conditional tail index estimation}\label{parag2.3}

The important question of how to estimate $\rho_x$ from the multivariate
random sample ${\mathcal{X}}_n$ is
very similar to the problem of estimating the so-called \textit{extreme
value index,} which is based
on a sample of \textit{univariate} random variables. An attractive
estimation method has been proposed by \cite{Pick}, which can be
easily adapted
to our conditional approach: let $k=k_n$ be a sequence of integers
tending to infinity and let $k/n \to0$ as $n \to\infty$. A
Pickands-type estimate of $\rho_x$ can be derived as
\begin{eqnarray*}
\hat\rho_x=\log2 \biggl(\log\frac
{\hat{\varphi}_{1-{(2k-1)}/{(n\hat{F}_{X} (x))}}(x)- \hat{\varphi
}_{1-{(4k-1)}/{(n\hat{F}_{X} (x))}}(x)}
{\hat{\varphi}_{1-{(k-1)}/{(n\hat{F}_{X} (x))}}(x)- \hat{\varphi
}_{1-{(2k-1)}/{(n\hat{F}_{X} (x))}}(x)}\biggr)^{-1} .
\end{eqnarray*}
The following result is particularly important since it allows the
hypothesis $\rho_x > 0$ to be tested and will later be employed to derive
asymptotic confidence intervals for $\varphi(x)$.

\begin{Theorem}
\label{thm2}
\textup{(i)} If \eqref{(2.2)} holds, $k_n\to\infty$ and $k_n/n \to
0$, then $\hat\rho_x \stackrel{p}{\longrightarrow} \rho_x$.
{\smallskipamount=0pt
\begin{enumerate}[(iii)]
\item[(ii)] If \eqref{(2.2)} holds, $k_n/n \to0$ and $k_n/\log\log
n\to\infty$, then $\hat\rho_x \stackrel{\mathit
{a.s.}}{\longrightarrow} \rho_x$.
\item[(iii)] Assume that $U(t):=\varphi_{1-{1}/{(tF_X(x))}}(x)$,
$t>\frac{1}{F_X(x)}$, has a positive derivative and that there exists
a positive
function $A(\cdot)$ such that for $z>0$,
$
\lim_{t\to\infty} \{(tz)^{1+1/\rho_x}  U'(tz) - t^{1+1/\rho
_x} U'(t)\}/A(t) = \pm\log(z),
$
for either choice of the sign \textup{(}$\Pi$-variation, which will in the sequel
be denoted by: $\pm
t^{1+1/\rho_x}U'(t) \in\Pi(A)$\textup{)}. Then,
%
\begin{eqnarray}
\label{Nrho1}
\sqrt{k_n} (\hat\rho_x - \rho_x) \stackrel{{d}}{\longrightarrow}
{\mathcal{N}
}(0,\sigma^2(\rho_x)) ,
\end{eqnarray}
with asymptotic variance $\sigma^2(\rho_x)=\rho^{2}_x(2^{1-{2}/{\rho_x}} + 1)/\{(2^{-{1}/{\rho_x}}-1)\log4\}^2$, for
$k_n\to\infty$ satisfying $k_n = \mathrm{o}(n/g^{-1}(n))$, where $g^{-1}$ is
the generalized inverse function of $g(t)= t^{3+ {2}/{\rho
_x}}\{ U'(t)/A(t)\}^2$.
\item[(iv)] If, for some $\kappa>0$ and $\delta>0$, the function
$\{ t^{\rho_x -1}F'(\varphi(x)-\frac{1}{t} \mid  x) - \delta
\} \in
 \mathit{RV} _{-\kappa}$, then \eqref{Nrho1} holds with $g(t)= t^{3+
{2}/{\rho_x}}\{ U'(t)/(t^{1+{1}/{\rho
_x}}U'(t)-[\delta F_X(x)]^{-1/\rho_x}(\rho_x)^{{1}/{\rho
_x}-1})\}^2$.
\end{enumerate}}
\end{Theorem}

\begin{Remark}
Note that the second order regular variation conditions (iii) and
(iv) of
Theorem \ref{thm2} are difficult to check in practice, which makes the
theoretical choice of the sequence $\{k_n\}$ a hard problem. In
practice, in order
to choose a reasonable estimate $\hat\rho_x(k_n)$ of $\rho_x$, one
can construct
the plot of $\hat\rho_x$, consisting of the points $\{(k,\hat\rho
_x(k)),1\leq
k<n\hat{F}_X(x)/4\}$, and select a value of $\rho_x$ at which the obtained
graph looks stable. This technique is known as the \textit{Pickands
plot} in the
univariate extreme value literature (see, for example, \cite{RES2} and
the references therein, Section 4.5,
pages 93--96).
This is this kind of idea which guides the
automatic data-driven rule we suggest in Section \ref{sec3}.
\end{Remark}

We can also easily adapt the well-known moment estimator for the index
of a
univariate extreme value distribution (Dekkers \textit{et al.}~\cite
{DEKK}) to our
conditional setup. Define
\begin{eqnarray*}
&&M^{(j)}_n=\frac{1}{k}\sum_{i=0}^{k-1}\bigl(\log\hat{\varphi
}_{1-{i}/{(n\hat{F}_{X}(x))}}(x)- \log\hat{\varphi}_{1-{k}/{(n\hat{F}_{X}
(x))}}(x)\bigr)^j\\
&&\quad  \mbox{for each }  j=1,2 \mbox{ and }  k=k_n<n .
\end{eqnarray*}
We can then define the moment-type estimator for the conditional
regular-variation exponent $\rho_x$ as
\[
\tilde\rho_x=-\biggl\{ M^{(1)}_n +1 -\frac{1}{2}\bigl[1-
\bigl(M^{(1)}_n\bigr)^2/M^{(2)}_n\bigr]^{-1}\biggr\}^{-1} .
\]

\begin{Theorem}
\label{thm?}
\textup{(i)} If \eqref{(2.2)} holds, $k_n/n \to0$ and $k_n\to\infty
$, then $\tilde\rho_x \stackrel{p}{\longrightarrow} \rho_x$.
{\smallskipamount=0pt
\begin{enumerate}[(iii)]
\item[(ii)] If \eqref{(2.2)} holds, $k_n/n \to0$ and $k_n/(\log
n)^{\delta}\to\infty$ for some $\delta>0$, then $\tilde\rho_x
\stackrel{\mathit{a.s.}}{\longrightarrow} \rho_x$.

\item[(iii)] If $\pm t^{1/\rho_x}\{\varphi(x)-U(t)\}\in\Pi(B)$ for
some positive function $B$, then
$
\sqrt{k_n} (\tilde\rho_x - \rho_x)
$
has, asymptotically, a normal distribution with mean zero and variance
\[
\rho_x(2+\rho_x)(1+\rho_x)^2\biggl\{4-8\frac{(2+\rho_x)}{(3+\rho
_x)}+\frac{(11+5\rho_x)(2+\rho_x)}{(3+\rho_x)(4+\rho_x)}\biggr\}
\]
for $k_n\to\infty$ satisfying $k_n = \mathrm{o}(n/g^{-1}(n))$, where
$g(t)= t^{1+{2}/{\rho_x}}[\{\log\varphi(x)-\break \log U(t)\}
/B(t)]^2$.
\end{enumerate}}
\end{Theorem}

\begin{Remark}
Note that the $\Pi$-variation condition $\pm t^{1+{1}/{\rho
_x}}U'(t)\in\Pi$
of Theorem \ref{thm2}(iii) is equivalent to $\pm(t^{1/\rho_x}\{
\varphi(x)-U(t)\})' \in \mathit{RV} _{-1}$,
following Theorem A.3 in \cite{DDH}, and that this equivalent regular-variation
condition implies that $\pm t^{1/\rho_x}\{\varphi(x)-U(t)\} \in\Pi
$, according to
\cite{RES}, Proposition 0.11(a), with auxiliary function $B(t)=\pm t
(t^{1/\rho_x}\{\varphi(x)-U(t)\})'$. Hence, the condition of Theorem
\ref{thm2}(iii) implies that of Theorem \ref{thm?}(iii).
Note, also, that a result similar to Theorem~\ref{thm?}(iii) can be
stated under
the conditions of Theorem \ref{thm2}(iv).
\end{Remark}

\subsection{Asymptotic confidence intervals}

The next theorem enables the construction of confidence intervals for
$\varphi(x)$ and for high
quantile-type frontiers $\varphi_{1-p_n/F_X(x)}(x)$ when $p_n\to0$
and $np_n\to\infty$.

\begin{Theorem}
\label{thm3}
{\smallskipamount=0pt
\begin{enumerate}[(iii)]
\item[(i)] Suppose that $F(\cdot|x)$ has a positive density $F'(\cdot
|x)$ such that $F'(\varphi(x)-\frac{1}{t}\mid x) \in \mathit{RV} _{1-\rho
_x}$. Then,
\begin{eqnarray*}
\sqrt{2k_n} \frac{\hat\varphi_{1-{(k_n-1)}/{(n\hat F_X(x))}}(x)-
\varphi_{1-{p_n}/{F_X(x)}}(x) }{\hat\varphi_{1-
{(k_n-1)}/{(n\hat F_X(x))}}(x)- \hat\varphi_{1-{(2k_n-1)}/{(n\hat
F_X(x))}}(x) }
\stackrel{{d}}{\longrightarrow} {\mathcal{N}}(0,V_1(\rho_x)) ,
\end{eqnarray*}
where $V_1(\rho_x) = \rho_x^{-2} 2^{1-2/\rho_x}/(2^{-1/\rho
_x}-1)^2$, provided that $p_n \to0$, $np_n\to\infty$ and $k_n=[np_n]$.

\item[(ii)] Suppose that the conditions of Theorem \textup{\ref{thm2}(iii)} or
\textup{(iv)} hold, and define
\begin{eqnarray*}
\hat\varphi^*_1(x) &:=& (2^{1/\hat\rho_x} - 1)^{-1}
\bigl\{\hat\varphi_{1-{(k_n-1)}/{(n\hat F_X(x))}}(x) - \hat\varphi
_{1-{(2k_n-1)}/{(n\hat F_X(x))}}(x)\bigr\}\\
&&{}+ \hat\varphi_{1-{(k_n-1)}/{(n\hat F_X(x))}}(x) .
\end{eqnarray*}
Then, putting $V_2(\rho_x) = 3\rho_x^{-2}2^{-1-2/\rho_x}/(2^{-1/\rho
_x}-1)^6$,
we have
\begin{eqnarray*}
\sqrt{2k_n} \frac{\hat\varphi^*_1(x) - \varphi(x) }
{\hat\varphi_{1-{(k_n-1)}/{(n\hat F_X(x))}}(x)- \hat\varphi
_{1-{(2k_n-1)}/{(n\hat F_X(x))}}(x) }
\stackrel{{d}}{\longrightarrow} {\mathcal{N}}(0,V_2(\rho_x)) .
\end{eqnarray*}
\item[(iii)] Suppose that the conditions of Theorem \textup{\ref{thm2}(iii)}
or \textup{(iv)} hold, and define
\begin{eqnarray*}
\tilde\varphi^*_1(x) &:=& (2^{1/\rho_x} - 1)^{-1}\bigl\{
\hat\varphi_{1-{(k_n-1)}/{(n\hat F_X(x))}}(x) - \hat\varphi
_{1-{(2k_n-1)}/{(n\hat F_X(x))}}(x)\bigr\}\\
&&{}+ \hat\varphi_{1-{(k_n-1)}/{(n\hat F_X(x))}}(x) .
\end{eqnarray*}
Then, putting $V_3(\rho_x) = \rho_x^{-2}2^{-2/\rho_x}/(2^{-1/\rho
_x}-1)^4$, we have
\begin{eqnarray}\label{2.12.bis}
&&\sqrt{2k_n} \frac{\tilde\varphi^*_1(x) - \varphi(x) }
{\hat\varphi_{1-{(k_n-1)}/{(n\hat F_X(x))}}(x)- \hat\varphi
_{1-{(2k_n-1)}/{(n\hat F_X(x))}}(x) }\nonumber\\
&&\quad \stackrel{{d}}{\longrightarrow} {\mathcal{N}}(0,V_3(\rho_x))
,\nonumber\\[-8pt]\\[-8pt]
&&\bigl\{\hat\varphi_{1-{(k_n-1)}/{(n\hat F_X(x))}}(x)-
\hat\varphi_{1-{(2k_n-1)}/{(n\hat
F_X(x))}}(x)\bigr\}\Big/\biggl\{\frac{n}{2k_n}U'\biggl(\frac{n}{2k_n}\biggr)\biggr\}\nonumber\\
&&\quad \stackrel
{p}{\longrightarrow}\rho_x(1-2^{-1/\rho_x}) .\nonumber
\end{eqnarray}
\end{enumerate}}
\end{Theorem}

\begin{Remark}
\label{rmk424}
Note that Theorem \ref{thm3}(ii) is still valid if the
estimate $\hat\rho_x$ is replaced by the true value $\rho_x$, up to
a change of
the asymptotic variance.
It is easy to see
that $V_2(\rho_x)\geq V_3(\rho_x)$ and so the estimator $\tilde
\varphi^*_1(x)$
of $\varphi(x)$ is asymptotically more efficient than $\hat\varphi^*_1(x)$.
We also conclude from \eqref{2.12.bis} that $\tilde\varphi^*_1(x)$ and
$\hat\varphi^*_1(x)$ have the same rate of convergence, namely
$nU'(\frac{n}{2k_n})/(2k_n)^{3/2}$. In the particular case where
$L_x(\{\varphi(x)-y\}^{-1})=\ell_x$ in \eqref{fxy}, we have
$U'(\frac{n}{2k_n})=\frac{1}{\rho_x}(\frac{1}{\ell_x})^{1/\rho
_x}(\frac{2k_n}{n})^{1+1/\rho_x}$.
Note, also, that in this particular case, the condition of Theorem~\ref
{thm3}(i) holds, that is, $F'(\varphi(x)-\frac{1}{t}\mid x)=
\frac{\ell_x \rho_x}{F_X(x)} (\frac{1}{t})^{\rho_x -
1} \in
 \mathit{RV} _{1-\rho_x}$. However, the conditions of Theorem~\ref{thm2}(iii)
and (iv) do not hold\vspace*{-1pt} since both functions $t^{1+ {1}/{\rho_x}}U'(t)=
\frac{1}{\rho_x}(\frac{1}{\ell_x})^{1/\rho_x}$ and
$t^{{\rho_x}-1}F'(\varphi(x) - \frac{1}{t}\mid x)=\frac{\ell_x\rho
_x}{F_X(x)}$
are constant\vspace*{1pt} in $t$. Nevertheless,
the conclusions of Theorem~\ref{thm2}(iii) and (iv) hold in this case
for all sequences $k_n \to\infty$ satisfying $\frac{k_n}{n} \to
0$. The same is true for the conclusion of Theorem~\ref{thm3}(ii).
\end{Remark}

\begin{Theorem}
\label{thm3bis}
If the condition of Corollary \ref{(corcor)} holds, $k_n \to\infty$ and
$k_n/n\to0$ as $n\to\infty$, then
\begin{eqnarray*}
&&\{\rho_xk^{1/2}_n/(k_n/n\ell_x)^{1/\rho_x}
\} \bigl[\hat{\varphi}_{1-(k_n-1)/(n\hat{F}_{X}(x))}(x) +
(k_n/n\ell_x)^{1/\rho_x} - \varphi
(x)\bigr]\\
&&\quad {} \stackrel{d}{\longrightarrow}
\mathcal{N}(0,1)\qquad  \mbox{as }  n \rightarrow\infty .
\end{eqnarray*}
\end{Theorem}

\begin{Remark}\label{rmq29}
The optimization of the asymptotic mean-squared error of\break
$\hat{\varphi}_{1-(k_n-1)/(n\hat{F}_{X}(x))}(x)$ is not an appropriate criteria
for selecting the optimal $k_n$ since the resulting value of $k_n$ does
not depend on $n$.
\end{Remark}

We shall now construct asymptotic confidence intervals for both
$\varphi(x)$
and $\varphi_{1-p_n/F_X(x)}(x),$ using the sums $M^{(1)}_n$ and $M^{(2)}_n$.

\begin{Theorem}
\label{thm2.7}
\begin{enumerate}[(ii)]
\item[(i)] Under the conditions of Theorem \ref{thm3}\textup{(i)},
\begin{eqnarray*}
\sqrt{k_n} \frac{\hat\varphi_{1-{k_n}/{(n\hat F_X(x))}}(x)-
\varphi_{1-{p_n}/{F_X(x)}}(x) }{M^{(1)}_n\hat\varphi_{1-
{k_n}/{(n\hat F_X(x))}}(x)}
\stackrel{{d}}{\longrightarrow} {\mathcal{N}}(0,V_4(\rho_x)) ,
\end{eqnarray*}
where $V_4(\rho_x) = (1+1/\rho_x)^2$, provided that $p_n \to0$,
$np_n\to\infty$ and $k_n=[np_n]$.

\item[(ii)] Suppose that the conditions of Theorem \textup{\ref{thm?}(iii)}
hold and that
$U(\cdot)$ has a regularly varying derivative $U'\in\mathit{RV}_{-\rho_x}$.
Define the moment estimator
$\hat\varphi(x)=\hat\varphi_{1-k_n/(n\hat F_X(x))}(x)\{
1+M^{(1)}_n(1+\tilde\rho_x)\}.$
Then,
\begin{eqnarray*}
&\displaystyle\sqrt{k_n} \frac{\hat\varphi(x) - \varphi(x) }
{M^{(1)}_n(1+1/\tilde\rho_x)\hat\varphi_{1-k_n/(n\hat F_X(x))}(x)}
\stackrel{{d}}{\longrightarrow} {\mathcal{N}}(0,V_5(\rho_x)),&\\
&\hspace*{-10pt}\displaystyle V_5(\rho_x) =\rho^{2}_x\biggl[
\frac{\rho_x}{(2+\rho_x)}+\rho_x(2+\rho_x)\biggl\{4-8\frac
{(2+\rho_x)}{(3+\rho_x)}+\frac{(11+5\rho_x)(2+\rho_x)}{(3+\rho
_x)(4+\rho_x)}\biggr\}&\\
&\hspace*{-78pt}{}-\dfrac{4\rho_x}{(3+\rho_x)}\biggr] .\hspace*{78pt}&
\end{eqnarray*}
\end{enumerate}
\end{Theorem}
%

\subsection{Examples}
\begin{Example}
\label{exemple1}
We consider the case where the support frontier is
linear. We choose $(X,Y)$ uniformly distributed over the region
$D=\{(x,y)\mid  0\le x\le1,  0\le y \le x \}$. In this case (see, for
example, \cite{DNN}), it is easy to see
that $\varphi(x)=x$ and $F_X(x)[1-F(y\vert x)]=(\varphi(x)-y)^2$ for
all $0\le
y \le\varphi(x)$. Thus, $L_x(\cdot)=\ell_x=1$ and $\rho_x=2$ for
all $x$.
Therefore, the conclusions of all Theorems \ref{thm1}--\ref{thm3bis}
hold (see Remark \ref{rmk424}).
\end{Example}

\begin{Example}
\label{exemple2}
We now choose a nonlinear monotone upper boundary given by the Cobb--Douglas
model $Y=X^{1/2}\exp(-U)$, where $X$ is uniform on $[0,1]$ and $U$,
independent of~$X$, is exponential with parameter $\lambda=3$ (see,
for example, \cite{DNN}). Here, the frontier function is
$\varphi(x)=x^{1/2}$ and the conditional distribution function is
$F(y|x)=3x^{-1}y^2-2x^{-3/2}y^3$ for $0<x\leq1$ and $0\leq y\leq
\varphi(x)$. It is then easily seen that the extreme value condition
\eqref{(2.2)} or, equivalently, \eqref{fxy} holds with $\rho_x=2$ and
$L_x(z)=F_X(x)[3\varphi(x)-\frac{2}{z}]/[\varphi(x)]^3$ for all
$x\in\,]0,1]$
and $z>0$.
\end{Example}
%

\section{Finite-sample performance}\label{sec3}

The simulation experiments of this section illustrate how the convergence
results work in practice. We also apply our approach to a real data set on
the production activity of the French postal services.

\subsection{Monte Carlo experiment}

We will simulate 2000 samples of size $n=5000$ according the scenario
of Example \ref{exemple1} above. Here, $\varphi(x)=x$ and $\rho
_x=2$. Denote by $N_x=n\hat F_X(x)$ the number of observations
$(X_i,Y_i)$ with $X_i \le x$. By construction of the estimators $\hat
\rho_x$ and $\hat\varphi^*_1(x)$, the threshold $k_n(x)$ can vary
between 1 and $N_x/4$. For the estimator with known $\rho_x$ and
$\tilde\varphi^*_1(x)$, $k_n(x)$ is bounded by $N_x/2$ and, finally,
for the moment estimators $\tilde\rho_x$ and $\hat\varphi(x)$, the
upper bound for $k_n(x)$ is given by $N_x-1$.
So, in our Monte Carlo experiments for the Pickands estimator, $k_n(x)$
was selected on a grid of values determined by the observed value of
$N_x$. We choose $k_n(x)=[N_x/4] - k + 1$, where $k$ is an integer
varying between 1 and $[N_x/4]$. In the tables below, $\bar N_x$ is
the average value observed over the 2000 Monte Carlo replications. The
tables display the values of $\bar{k}_n(x)$, which is the average of
the Monte Carlo values of $k_n(x)$ obtained for a fixed selection of
values of $k$. For the moment estimators, the upper values of $k_n(x)$
were chosen as $N_x-1$. The tables display only a part of the results
to save space, but in each case, we typically choose a set of values of
$k$ that includes not only the most favorable cases, but also covers a
wide range of values for $k_n(x)$.
These tables
provide the Monte Carlo estimates of the bias and the mean-squared
error (MSE)
of the various estimators computed over the 2000 random replications,
as well
as the average lengths and the achieved coverages of the corresponding
95\%
asymptotic confidence intervals. They display only the results for
$x$ ranging over $\{0.25,0.5,1\},$ to save space.

\begin{table}
\caption{Pickands and known $\rho_x$ cases: bias ($B$) and
mean-squared error ($\mathit{MSE}$) of the estimates}
\label{tab:MCex2_BMSE}
\begin{tabular*}{\textwidth}{@{\extracolsep{\fill}}d{4.1}d{2.5}d{5.5}d{2.5}d{2.5}d{2.5}d{1.5}@{}}
\hline
\multicolumn{1}{@{}l}{$\bar{k}_n(x)$}&\multicolumn{1}{l}{$B_{\hat\rho_x}$}&\multicolumn{1}{l}{$\mathit{MSE}_{\hat\rho_x}$}&\multicolumn{1}{l}{$B_{\hat
\varphi^*_1(x)}$}&\multicolumn{1}{l}{$\mathit{MSE}_{\hat\varphi^*_1(x)}$}&\multicolumn{1}{l}{$B_{\tilde
\varphi^*_1(x)}$}&\multicolumn{1}{l@{}}{$\mathit{MSE}_{\tilde\varphi^*_1(x)}$} \\
\hline
\multicolumn{7}{@{}l}{$x= 0.25,\bar N_x=312$, FDH: $B_{\hat\varphi_1(x)}=-0.012591$, $\mathit{MSE}_{\hat\varphi_1(x)}=0.000203$ }\\
[3pt]
77.7& -0.25757 & 784.19539 & -0.02585 & 6.93961 & 0.00021 & 0.00028 \\
74.4& 0.41215 & 17.20703 & 0.03723 & 0.14471 & 0.00024 & 0.00028 \\
71.0& 0.42344 & 105.75775 & 0.03830 & 0.89895 & 0.00016 & 0.00028 \\
67.7& 0.44401 & 16.30552 & 0.03877 & 0.11468 & 0.00030 & 0.00028 \\
64.4& 0.30552 & 145.08207 & 0.02564 & 1.01166 & 0.00031 & 0.00029 \\
61.0& 0.68905 & 35.13730 & 0.05654 & 0.24012 & 0.00053 & 0.00029 \\
57.7& 0.82177 & 15489.98302 & 0.05929 & 89.02353 & 0.00053 & 0.00029 \\
54.3& 1.17914 & 1780.66037 & 0.08527 & 9.90370 & 0.00055 & 0.00029 \\
51.0& -4.41384 & 13169.38480 & -0.33207 & 74.80129 & 0.00046 & 0.00030\\
47.6& 0.03147 & 3204.61688 & -0.00179 & 14.27123 & 0.00064 & 0.00029 \\
[6pt]
\multicolumn{7}{@{}l}{$x= 0.50$, $\bar N_x=1250$, FDH: $B_{\hat\varphi_1(x)}=-0.012563$, $\mathit{MSE}_{\hat\varphi_1(x)}=0.000200$}\\
[3pt]
312.1& 0.09248 & 0.22503 & 0.01696 & 0.00735 & 0.00026 & 0.00029 \\
297.0& 0.09311 & 0.24340 & 0.01668 & 0.00759 & 0.00012 & 0.00029 \\
281.9& 0.09124 & 0.24958 & 0.01595 & 0.00742 & -0.00001 & 0.00029 \\
266.8& 0.09201 & 0.27538 & 0.01579 & 0.00780 & -0.00009 & 0.00029 \\
251.7& 0.08954 & 0.29784 & 0.01490 & 0.00797 & -0.00042 & 0.00030 \\
236.6& 0.09840 & 0.33195 & 0.01584 & 0.00831 & -0.00049 & 0.00030 \\
221.5& 0.11387 & 0.38048 & 0.01768 & 0.00893 & -0.00043 & 0.00030 \\
206.3& 0.12297 & 0.47557 & 0.01840 & 0.01038 & -0.00060 & 0.00030 \\
191.2& 0.12060 & 0.43562 & 0.01720 & 0.00881 & -0.00081 & 0.00030 \\
176.1& 0.14573 & 0.72946 & 0.01989 & 0.01371 & -0.00080 & 0.00029 \\
[6pt]
\multicolumn{7}{@{}l}{$x= 1.00$, $\bar N_x=5000$, FDH: $B_{\hat\varphi_1(x)}=-0.012663$, $\mathit{MSE}_{\hat\varphi_1(x)}=0.000202$}\\
[3pt]
1250.0& 0.02755 & 0.04085 & 0.01025 & 0.00540 & 0.00078 & 0.00028 \\
1188.0& 0.02863 & 0.04254 & 0.01047 & 0.00537 & 0.00085 & 0.00028 \\
1126.0& 0.02780 & 0.04643 & 0.00991 & 0.00557 & 0.00065 & 0.00029 \\
1064.0& 0.02689 & 0.05068 & 0.00953 & 0.00575 & 0.00064 & 0.00030 \\
1002.0& 0.02890 & 0.05241 & 0.00981 & 0.00559 & 0.00061 & 0.00029 \\
940.0& 0.02670 & 0.05545 & 0.00875 & 0.00552 & 0.00032 & 0.00029 \\
878.0& 0.02738 & 0.06064 & 0.00872 & 0.00564 & 0.00029 & 0.00029 \\
816.0& 0.02877 & 0.06738 & 0.00882 & 0.00577 & 0.00024 & 0.00028 \\
754.0& 0.03001 & 0.07071 & 0.00899 & 0.00562 & 0.00037 & 0.00028 \\
692.0& 0.03686 & 0.07869 & 0.01065 & 0.00583 & 0.00065 & 0.00029 \\
\hline
\end{tabular*}
\end{table}

We will first comment on the results obtained for the Pickands
estimators and for
the estimator of $\varphi(x)$ obtained with the knowledge that $\rho
_x=p+1=2$ (the jump of
the joint density of $(X,Y)$ at the frontier); these results are
displayed in Tables \ref{tab:MCex2_BMSE} and
\ref{tab:MCex2_CI}. We observe that the Pickands estimates $\hat\rho
_x$ and
$\hat\varphi^*_1(x)$ behave much better when the sample size $N_x$ increases,
although the convergence is rather slow. In contrast, even with the
smallest sample size $N_x$
(for $x=0.25$), the estimator $\tilde\varphi^*_1(x)$ computed with
the true
value of $\rho_x=2$ provides remarkable estimates of $\varphi(x)$ and
is rather stable with respect to the choice of $k_n(x)$. We see the
improvement of $\tilde\varphi^*_1(x)$ over the FDH in terms of the
bias, without significantly increasing the MSE. The achieved coverages
of the normal confidence intervals obtained from $\tilde\varphi
^*_1(x)$ are also quite satisfactory and much easier to derive than
those obtained from the FDH estimator.
As soon as $N_x$ is greater than 1000, all of the estimators provide
reasonably good confidence intervals of the corresponding unknown, with
quite good achieved coverages. In these cases ($N_x \ge1000$), we also
observe some stability of the results with respect to the choice of $k_n(x)$.

\begin{table}
\caption{Pickands and known $\rho_x$ cases: average lengths
($\mathit{avl}$) and coverages ($\mathit{cov}$) of the 95\% confidence intervals}
\label{tab:MCex2_CI}
\begin{tabular*}{\textwidth}{@{\extracolsep{\fill}}d{4.1}d{5.4}d{1.4}d{4.4}d{1.4}d{1.4}d{1.4}@{}}
\hline
\multicolumn{1}{@{}l}{$\bar{k}_n(x)$}&\multicolumn{1}{l}{$\mathit{avl}_{\hat\rho_x}$}&\multicolumn{1}{l}{$\mathit{cov}_{\hat\rho
_x}$}&\multicolumn{1}{l}{$\mathit{avl}_{\hat\varphi^*_1(x)}$}&\multicolumn{1}{l}{$\mathit{cov}_{\hat\varphi
^*_1(x)}$}&\multicolumn{1}{l}{$\mathit{avl}_{\tilde\varphi^*_1(x)}$}&\multicolumn{1}{l@{}}{$\mathit{cov}_{\tilde\varphi
^*_1(x)}$} \\
\hline
\multicolumn{7}{@{}l}{$x= 0.25$, $\bar N_x=312$}\\
[3pt]
77.7& 630.9019 & 0.9040 & 59.3041 & 0.8925 & 0.0670 & 0.9455 \\
74.4& 18.4635 & 0.9060 & 1.6821 & 0.8970 & 0.0670 & 0.9505 \\
71.0& 92.5814 & 0.9000 & 8.5104 & 0.8960 & 0.0670 & 0.9480 \\
67.7& 18.6125 & 0.8990 & 1.5673 & 0.8910 & 0.0670 & 0.9485 \\
64.4& 131.0169 & 0.8910 & 10.9372 & 0.8845 & 0.0670 & 0.9525 \\
61.0& 37.9315 & 0.8960 & 3.1260 & 0.8840 & 0.0671 & 0.9465 \\
57.7& 14491.7449 & 0.8965 & 1098.2578 & 0.8850 & 0.0671 & 0.9470 \\
54.3& 1735.9675 & 0.8930 & 129.3070 & 0.8820 & 0.0671 & 0.9430 \\
51.0& 13077.3352 & 0.8910 & 981.3170 & 0.8805 & 0.0671 & 0.9440 \\
47.6& 3374.6016 & 0.8925 & 224.7041 & 0.8735 & 0.0672 & 0.9410 \\
[6pt]
\multicolumn{7}{@{}l}{$x= 0.50$, $\bar N_x= 1250$}\\
[3pt]
312.1& 1.7798 & 0.9295 & 0.3232 & 0.9195 & 0.0670 & 0.9485 \\
297.0& 1.8330 & 0.9255 & 0.3248 & 0.9245 & 0.0669 & 0.9490 \\
281.9& 1.8810 & 0.9250 & 0.3247 & 0.9240 & 0.0669 & 0.9475 \\
266.8& 1.9457 & 0.9220 & 0.3269 & 0.9240 & 0.0669 & 0.9460 \\
251.7& 2.0095 & 0.9200 & 0.3279 & 0.9145 & 0.0668 & 0.9505 \\
236.6& 2.1038 & 0.9195 & 0.3329 & 0.9165 & 0.0668 & 0.9420 \\
221.5& 2.2256 & 0.9150 & 0.3409 & 0.9100 & 0.0668 & 0.9390 \\
206.3& 2.3707 & 0.9115 & 0.3506 & 0.9075 & 0.0668 & 0.9440 \\
191.2& 2.4375 & 0.9105 & 0.3468 & 0.9085 & 0.0667 & 0.9455 \\
176.1& 2.7460 & 0.9155 & 0.3754 & 0.9080 & 0.0667 & 0.9440 \\
[6pt]
\multicolumn{7}{@{}l}{$x= 1.00$, $\bar N_x=5000$}\\
[3pt]
1250.0& 0.8019 & 0.9645 & 0.2909 & 0.9605 & 0.0670 & 0.9540 \\
1188.0& 0.8238 & 0.9625 & 0.2914 & 0.9595 & 0.0670 & 0.9555 \\
1126.0& 0.8463 & 0.9535 & 0.2914 & 0.9495 & 0.0670 & 0.9425 \\
1064.0& 0.8707 & 0.9510 & 0.2915 & 0.9445 & 0.0670 & 0.9435 \\
1002.0& 0.8994 & 0.9530 & 0.2922 & 0.9455 & 0.0670 & 0.9475 \\
940.0& 0.9273 & 0.9445 & 0.2918 & 0.9420 & 0.0669 & 0.9460 \\
878.0& 0.9614 & 0.9420 & 0.2923 & 0.9450 & 0.0669 & 0.9420 \\
816.0& 1.0002 & 0.9450 & 0.2932 & 0.9440 & 0.0669 & 0.9500 \\
754.0& 1.0426 & 0.9475 & 0.2939 & 0.9460 & 0.0669 & 0.9550 \\
692.0& 1.0976 & 0.9455 & 0.2966 & 0.9430 & 0.0670 & 0.9455 \\
\hline
\end{tabular*}
\end{table}

\begin{table}
\tabcolsep=0pt
\caption{Moment Estimators: bias, MSE, average lengths and coverages}
\label{tab:MCex2_mom}
\begin{tabular*}{\textwidth}{@{\extracolsep{\fill}}d{4.1}d{2.5}d{4.5}d{2.5}d{1.5}d{5.4}d{1.4}d{3.4}d{1.4}@{}}
\hline
\multicolumn{1}{@{}l}{$\bar{k}_n(x)$}&\multicolumn{1}{l}{$B_{\tilde\rho_x}$}&\multicolumn{1}{l}{$\mathit{MSE}_{\tilde\rho
_x}$}&\multicolumn{1}{l}{$B_{\hat\varphi(x)}$}&\multicolumn{1}{l}{$\mathit{MSE}_{\hat\varphi
(x)}$}&\multicolumn{1}{l}{$\mathit{avl}_{\tilde\rho_x}$}&\multicolumn{1}{l}{$\mathit{cov}_{\tilde\rho_x}$}&\multicolumn{1}{l}{$\mathit{avl}_{\hat
\varphi(x)}$}&\multicolumn{1}{l@{}}{$\mathit{cov}_{\hat\varphi(x)}$} \\
\hline
\multicolumn{9}{@{}l}{$x= 0.25$, $\bar N_x=312$}\\
[3pt]
150.4 & 0.36520 & 1.47278 & -0.04187 & 0.00339 & 2.5969 & 0.8900 &0.0869 & 0.3350 \\
137.9 & 0.35077 & 1.86333 & -0.03615 & 0.00337 & 2.8243 & 0.8905 &0.0939 & 0.3765 \\
125.3 & 0.33799 & 1.26492 & -0.03080 & 0.00226 & 2.7378 & 0.8990 &0.0893 & 0.4435 \\
112.9 & 0.30315 & 1.02334 & -0.02670 & 0.00173 & 2.7495 & 0.9005 &0.0874 & 0.4840 \\
100.4 & 0.27374 & 0.93872 & -0.02284 & 0.00139 & 2.8414 & 0.8930 &0.0873 & 0.5495 \\
87.9 & 0.28569 & 1.22921 & -0.01810 & 0.00137 & 3.1695 & 0.8965 &0.0936 & 0.5860 \\
75.4 & 0.30500 & 9.96907 & -0.01330 & 0.00806 & 7.3693 & 0.8865 &0.2075 & 0.6340 \\
62.9 & 0.26381 & 29.37920 & -0.01097 & 0.02156 & 17.2434 & 0.8880 &0.4629 & 0.6740 \\
50.5 & 0.51850 & 18.67121 & -0.00130 & 0.01090 & 14.4349 & 0.8780 &0.3524 & 0.7020 \\
38.0 & 0.53418 & 21.11753 & 0.00124 & 0.00956 & 18.2022 & 0.8645 &0.3897 & 0.7225 \\
19.2 & 0.62323 & 267.28452 & 0.00481 & 0.06789 & 246.3768 & 0.8430 &3.8848 & 0.7525 \\
12.9 & -0.30491 & 1266.44113 & -0.00977 & 0.30730 & 1431.7282 & 0.8150& 22.2514 & 0.7315 \\
[6pt]
\multicolumn{9}{@{}l}{$x= 0.50$, $\bar N_x=1250$}\\
[3pt]
600.5 & 0.16644 & 0.16966 & -0.09657 & 0.01004 & 0.9860 & 0.8375 &0.0645 & 0.0575 \\
550.5 & 0.16412 & 0.16874 & -0.08407 & 0.00776 & 1.0281 & 0.8590 &0.0667 & 0.0890 \\
500.4 & 0.16750 & 0.17596 & -0.07212 & 0.00588 & 1.0818 & 0.8735 &0.0691 & 0.1360 \\
450.5 & 0.17133 & 0.18419 & -0.06106 & 0.00440 & 1.1442 & 0.8970 &0.0715 & 0.2155 \\
400.5 & 0.16370 & 0.19777 & -0.05158 & 0.00334 & 1.2099 & 0.9085 &0.0733 & 0.2945 \\
350.5 & 0.15716 & 0.20738 & -0.04270 & 0.00250 & 1.2897 & 0.9225 &0.0751 & 0.3815 \\
300.5 & 0.16437 & 0.23740 & -0.03370 & 0.00182 & 1.4051 & 0.9335 &0.0778 & 0.4775 \\
250.4 & 0.15151 & 0.25663 & -0.02649 & 0.00137 & 1.5307 & 0.9430 &0.0794 & 0.5650 \\
200.5 & 0.13915 & 0.28167 & -0.01987 & 0.00101 & 1.7031 & 0.9415 &0.0811 & 0.6475 \\
150.5 & 0.12971 & 0.36589 & -0.01373 & 0.00082 & 1.9765 & 0.9305 &0.0836 & 0.7180 \\
50.5 & 0.29865 & 6.19391 & 0.00098 & 0.00356 & 6.8895 & 0.8895 & 0.1734& 0.8000 \\
13.0 & -0.58590 & 9410.59672 & -0.01445 & 1.57034 & 10243.4270 & 0.8150& 131.6029 & 0.7550 \\
[6pt]
\multicolumn{9}{@{}l}{$x= 1.00$, $\bar N_x=5000$}\\
[3pt]
2000.0 & 0.13502 & 0.05141 & -0.14729 & 0.02230 & 0.5207 & 0.7685 &0.0664 & 0.0000 \\
1800.0 & 0.13019 & 0.05132 & -0.12609 & 0.01649 & 0.5471 & 0.8140 &0.0682 & 0.0025 \\
1600.0 & 0.12099 & 0.04935 & -0.10701 & 0.01202 & 0.5765 & 0.8455 &0.0697 & 0.0145 \\
1400.0 & 0.11212 & 0.05190 & -0.08930 & 0.00855 & 0.6129 & 0.8595 &0.0712 & 0.0455 \\
1200.0 & 0.10555 & 0.05445 & -0.07261 & 0.00584 & 0.6593 & 0.8965 &0.0727 & 0.1055 \\
1000.0 & 0.09393 & 0.05677 & -0.05771 & 0.00388 & 0.7168 & 0.9180 &0.0740 & 0.2325 \\
800.0 & 0.07446 & 0.05965 & -0.04469 & 0.00251 & 0.7911 & 0.9245 &0.0748 & 0.3680 \\
600.0 & 0.07713 & 0.07992 & -0.03069 & 0.00148 & 0.9179 & 0.9310 &0.0771 & 0.5615 \\
400.0 & 0.06905 & 0.10581 & -0.01877 & 0.00087 & 1.1221 & 0.9415 &0.0790 & 0.7255 \\
200.0 & 0.07559 & 0.20770 & -0.00744 & 0.00059 & 1.6176 & 0.9365 &0.0830 & 0.8375 \\
100.0 & 0.09821 & 0.49803 & -0.00225 & 0.00067 & 2.4204 & 0.9095 &0.0896 & 0.8465 \\
50.0 & 0.15884 & 1.20953 & 0.00051 & 0.00083 & 3.9082 & 0.8920 & 0.1034& 0.8420 \\
\hline
\end{tabular*}
\end{table}

We now turn to the performances of the moment estimators $\tilde\rho
_x$ and
$\hat\varphi(x)$. The results are displayed in Table \ref{tab:MCex2_mom}.
Note that we used the same seed in the Monte Carlo experiments as the one
used for the preceding tables. Compared with the Pickands estimators
$\hat\rho_x$ and $\hat\varphi^*_1(x)$,
we observe here much more reasonable results in terms of the bias and
MSE of
the estimators $\tilde\rho_x$ and $\hat\varphi(x)$. In addition,
when $N_x$ increases, the results are
much less sensitive to the choice of $k_n(x)$ than for the Pickands
estimators. We also observe that the most favorable values of $k_n(x)$ for
estimating $\rho_x$ and $\varphi(x)$ are not necessarily in the same
range of
values.
We note that the confidence intervals for $\rho_x$ achieve quite reasonable
coverage as soon as $N_x$ is greater than, say, 1000. However, the
results for
the confidence intervals of $\varphi(x)$ obtained from the moment estimator
$\hat\varphi(x)$ are very poor, even when $N_x$ is as large as 5000.
A more
detailed analysis of the Monte Carlo results allows us to conclude that this
comes from an under-evaluation of the asymptotic variance of $\hat
\varphi(x)$
given in Theorem \ref{thm2.7}. Indeed, in most of the cases, the Monte Carlo
standard deviation of $\hat\varphi(x)$ was larger than the asymptotic
theoretical expression by a factor of the order 2--5 when $N_x$
equalled $1250$, and by
a factor of the order 1.3--1.7 when it equalled $5000$. So, the poor
behavior seems to improve
slightly when $N_x$ increases, but at a very slow rate.

We could say that using the Pickands estimators $\hat\rho_x$ and
$\hat\varphi^*_1(x)$ is only reasonable in our setup when $N_x$
is larger than, say, 1000. These estimators are highly sensitive to the
choice of $k_n(x)$. The moment estimators $\tilde\rho_x$ and $\hat
\varphi(x)$ have a much better behavior in terms of bias and MSE, and
a greater stability with respect to the choice of $k_n(x)$, even for
moderate sample sizes. When $N_x$ is very large ($N_x=5000$), $\hat
\rho_x$ and $\hat\varphi^*_1(x)$ become more accurate than the
moment estimators.
On the other hand, the confidence intervals of $\rho_x$ constructed
from the
asymptotic distribution of $\hat\rho_x$ provide more satisfactory
results than those derived from the limit distribution of $\tilde\rho
_x$ for large values of $N_x$, say, $N_x\geq1000$.
For inference purposes on the frontier function itself, the estimate of the
asymptotic variance of the moment estimator $\hat\varphi(x)$ does not
provide reliable confidence
intervals, even for relatively large values of $N_x$.
In the latter case, it would be better to use the confidence
intervals obtained from the asymptotic distribution of
the Pickands estimator $\hat\varphi^*_1(x)$.

So, in terms of bias and MSE
computed over the 2000 random replications, as well as the average
lengths and
the achieved coverages of the $95\%$ asymptotic confidence intervals, the
moment estimators of $\rho_x$
and $\varphi(x)$ are sometimes preferable to the Pickands estimators and
sometimes not. It is difficult to imagine one procedure being
preferable in all
contexts. Hence, a sensible practice is not to restrict the frontier
analysis to
one procedure, but rather to check that both Pickands and moment estimators
point toward similar conclusions.
However, when $\rho_x$ is known, we have remarkable results for
$\tilde\varphi^*_1(x)$, even when $N_x$ is small, including
remarkable properties
of the resulting normal confidence intervals, with great stability with
respect to the choice of $k_n(x)$. Recall that in most situations described
thus far in the econometric literature on frontier analysis, this tail index
$\rho_x$ is supposed to be known and equal to $p+1$ (here, $\rho
_x=2$): this
corresponds to the common assumption that there is a jump of the joint
density of $(X,Y)$ at the frontier.

This might suggest the following strategy with a real data set. If
$\rho_x$ is known (typically equal to $p+1$ if the assumption of a
jump at the
frontier is reasonable), then we can use the estimator~$\tilde\varphi^*_1(x)$. If, on the other hand, $\rho_x$ is
unknown, we could
consider using the following two-step estimator: first, estimate $\rho
_x$ (the
moment estimator of $\rho_x$ seems the more appropriate, unless $N_x$
is large
enough) and, second, use the estimator $\tilde\varphi^*_1(x)$, as if
$\rho_x$ were known, by substituting the estimated value $\tilde\rho
_x$ or
$\hat\rho_x$ in place of $\rho_x$.
In a situation involving a real data set,
the best approach is not to favor the moment or
the Pickands estimator of $\rho_x$ in the first step, but to compute
$\tilde\varphi^*_1(x)$ by
substituting in each of them, in the hope that the two resulting values of
$\tilde\varphi^*_1(x)$ point toward similar conclusions.

It should be clear that the two-step estimator
$\tilde\varphi^*_1(x)$, obtained by
substituting in $\hat\rho_x$, does not necessarily coincide with the Pickands
estimator $\hat\varphi^*_1(x),$ which is, instead, obtained by a
simultaneous
estimation of $\rho_x$ and $\varphi(x)$. Indeed, in our
Monte Carlo exercise, we have observed that the most favorable values
of $k_n(x)$ for estimating
$\rho_x$ and $\varphi(x)$ are not necessarily in the same range of values.
Thus, nothing guarantees that the selected value $k_n(x)$
when computing $\hat\rho_x$ in the first step is the same as the one selected
when computing~$\hat\varphi^*_1(x)$. Of course, when $N_x$ is
very large, the
two values of $k_n(x)$ are expected to be similar, but the idea in the two-step
procedure is to use the asymptotic results of the more efficient estimator
$\tilde\varphi^*_1(x)$ and not those of $\hat\varphi^*_1(x)$.
In the next section, we suggest an ad hoc procedure for determining
appropriate values of $k_n(x)$ with a real data set.

\subsection{A data-driven method for selecting $k_n(x)$}

The question of selecting the optimal value of $k_n(x)$ is still an
open issue and is not addressed here. We will simply suggest an
empirical rule that turns out to give reasonable estimates of the
frontier in the simulated samples above.

First, we have observed in our Monte Carlo exercise that the optimal
value for selecting $k_n(x)$ when estimating the index $\rho_x$ is not
necessarily the same as the value for estimating $\varphi(x)$.
The~idea is thus to select first, for each $x$ (in a chosen grid of
values), a grid of values for $k_n(x)$ for estimating $\rho_x$. For
the Pickands estimator $\hat\rho_x$, we choose $k_n(x)=[N_x/4] - k +
1$, where $k$ is an integer varying between 1 and $[N_x/4]$, and for
the moment estimator $\tilde\rho_x$, we choose $k_n(x)=N_x - k$,
where $k$ is an integer varying between 1 and $N_x$.
We then evaluate the estimator $\hat\rho_x(k)$ (resp., $\tilde\rho
_x(k)$) and select the $k$ where the
variation of the results is the smallest. We achieve this by computing the
standard deviations of $\hat\rho_x(k)$ (resp., $\tilde\rho_x(k)$)
over a `window' of $2\times[\sqrt{N_x/4}]$ (resp., $2\times[\sqrt
{N_x}]$) successive values of $k$. The value of $k$ where this standard
deviation is minimal defines the value of $k_n(x)$.

We follow the same procedure for selecting a value for $k_n(x)$ for
estimating the
frontier $\varphi(x)$ itself. Here, in all of the cases, we choose a
grid of values for $k_n(x)$ given by $k= 1,\ldots, [\sqrt{N_x}]$ and select
the $k$ where the variation of the results is the smallest. To achieve this
here, we compute the standard deviations of $\tilde\varphi^*_1(x)$ (resp.,
$\hat\varphi^*_1(x)$ and $\hat\varphi(x)$) over a `window' of size
$2\times
\max(3, [\sqrt{N_x}/20])$ (this corresponds to having a window large
enough to
cover around 10\% of the possible values of $k$ in the selected range of
values for $k_n(x)$). From now on, we only present illustrations for
$\tilde\varphi^*_1(x)$ to save space.

\begin{figure}

\includegraphics{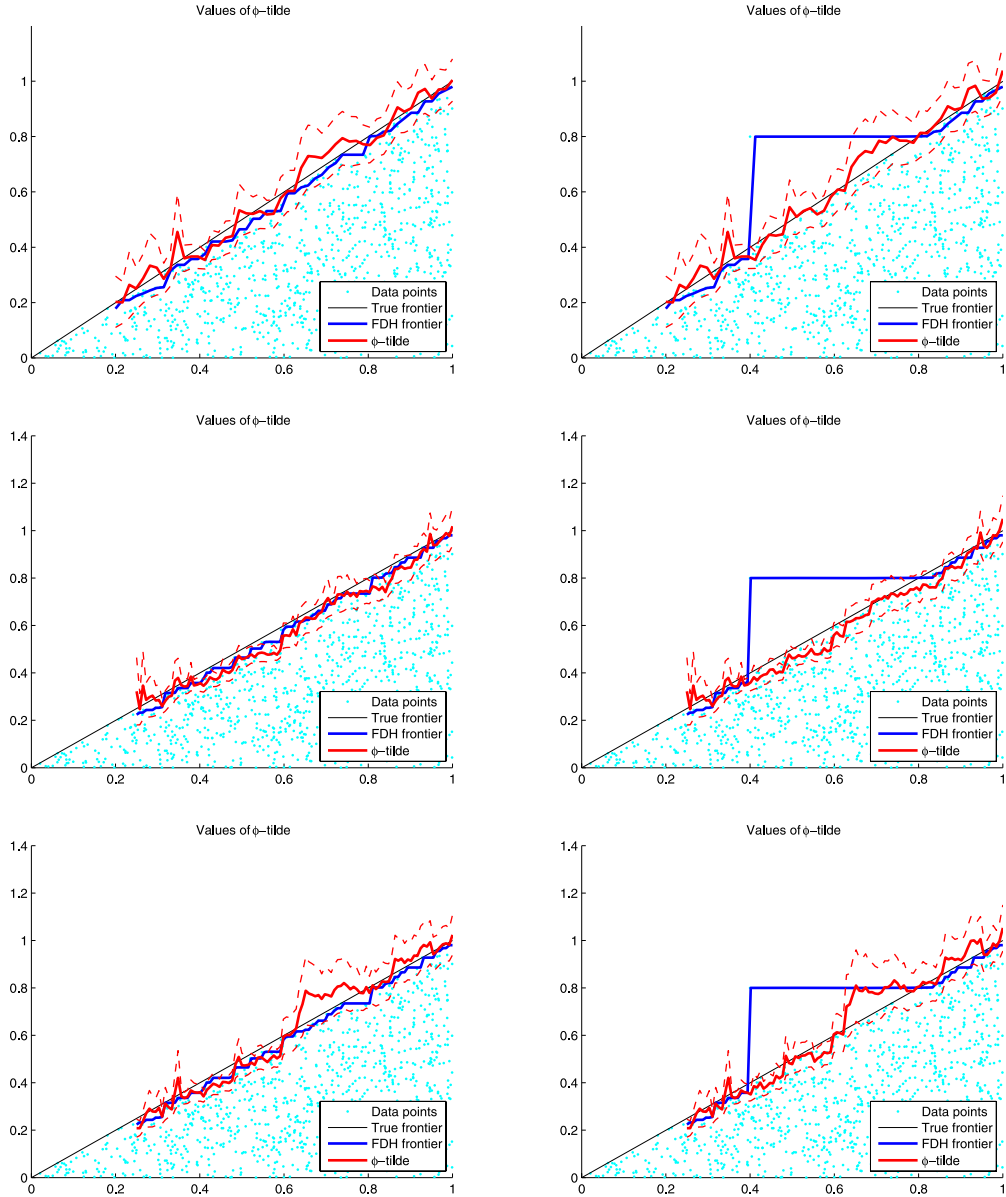}

\caption{Resulting estimator $\tilde\varphi^*_1(x)$ for a
uniform data set of size $n=1000$ (plus one outlier for the right
panels); from top to bottom, we have the cases $\rho_x=2$,
substituting in $\hat\rho_x$, substituting in $\tilde\rho_x$.}\label{fig1}
\end{figure}

For a sample generated with $n=1000$ in the uniform case, we get the
results shown in Fig.~\ref{fig1}.

In Fig. \ref{fig1}, the estimator $\tilde\varphi^*_1(x)$ is first computed
with the
true value $\rho_x=2$ (top panel of the figure), then with a plug-in value
of $\rho_x$ estimated by the Pickands estimator (middle panel) and
finally with a plug-in value
of $\tilde\rho_x$ estimated by the
moment estimator (bottom panel). The pointwise confidence
intervals are also displayed. The three right-hand panels correspond to
the same
data set plus one outlier. This allows us to see how our robust estimators
behave in the presence of outlying points, in contrast with the FDH estimator.
In particular, due to the remarkable behavior of $\tilde\varphi
^*_1(x)$ in the
Monte Carlo experiment, if we know that $\rho_x=2$, then we should use
the top
panel results and, according to our suggestion at the end of the preceding
section, if $\rho_x$ is unknown, we should use, in this particular
example, the bottom panel results, where
we replace $\rho_x$ by its moment estimator $\tilde\rho_x$ (since
here $N_x
\le1000$) and continue as if $\rho_x$ were known. It is quite
encouraging that
the two panels are very similar.

\begin{figure}[b]

\includegraphics{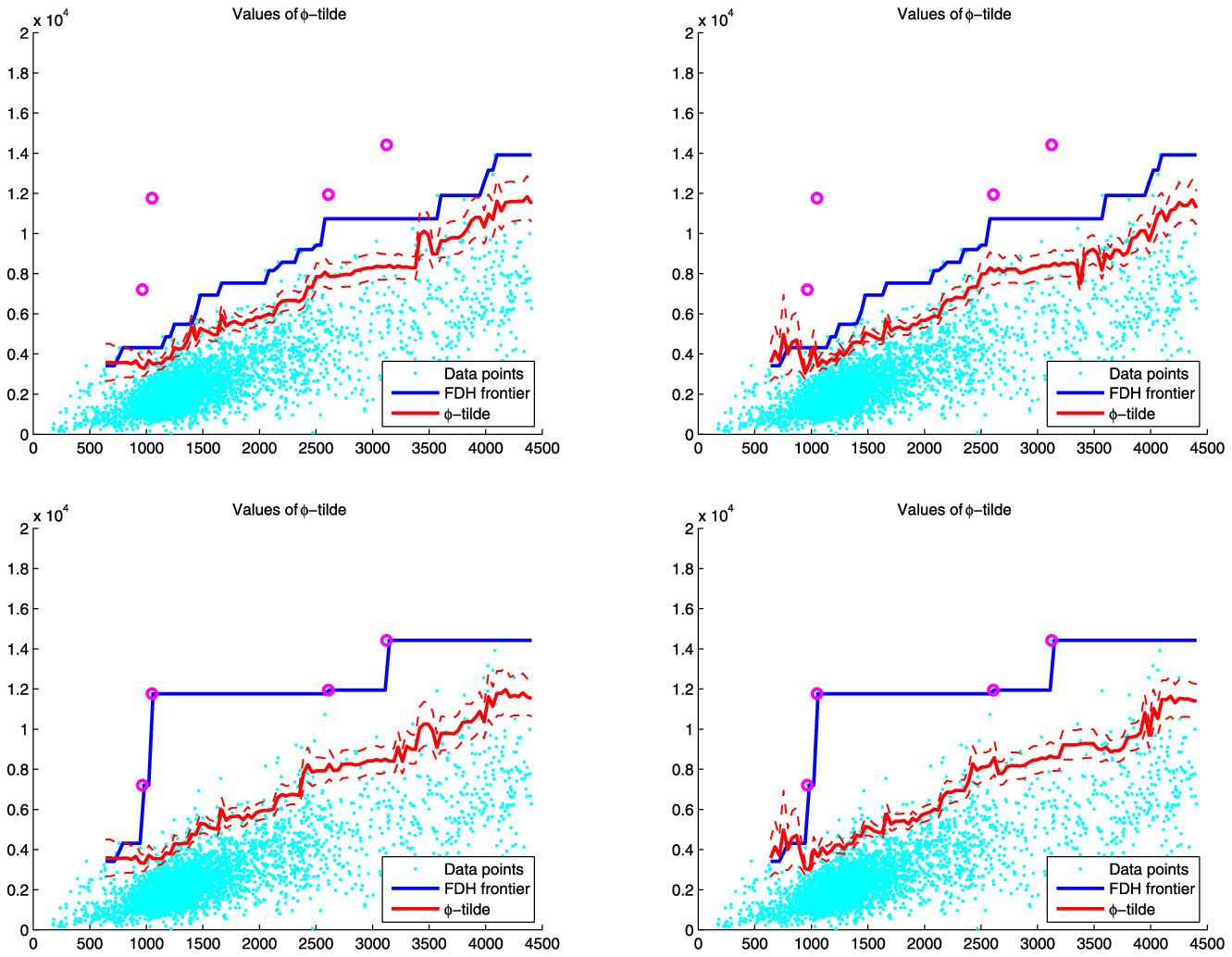}

\caption{The resulting estimator $\tilde\varphi^*_1(x)$ for the
French post offices. We include four extreme data points (circles) for
the bottom panels. From left to right, we have the cases $\rho_x=2$,
substituting in $\tilde\rho_x$.}\label{fig:poste1}
\end{figure}

\subsection{An application}

We use the same real data example as in \cite{CAZ}, which undertook
the frontier analysis of
$9521$ French post offices observed in $1994$, with $X$ as the quantity of
labor and $Y$ as the volume of delivered mail. In this illustration, we only
consider the $n=4000$ observed post offices with the smallest levels
$x_i$. We
used the empirical rules explained above for selecting reasonable
values for
$k_n(x)$. The cloud of points and the resulting estimates are provided in
Fig. \ref{fig:poste1}.

To save space, we only represent
$\tilde\varphi^*_1(x)$ when $\rho_x$ is supposed to be equal to 2 (left-hand
panels) and when it is estimated by the moment estimator (right-hand panels).
The FDH estimator is clearly determined by only a few very extreme
points. If we delete four extreme points from the sample (represented
by circles in the figure), then we obtain the pictures from the top
panels: the FDH estimator
changes drastically, whereas the extreme-value-based estimator
$\tilde\varphi^*_1(x)$ is very robust to the presence of these four extreme
points. We also note the considerable stability of the various forms of
the estimator
$\tilde\varphi^*_1(x)$.

\section{Concluding remarks}\label{sec4}

In our approach, we provide the necessary and sufficient
condition for the FDH estimator $\hat{\varphi}_1(x)$ to converge in
distribution, we specify
its asymptotic distribution with the appropriate convergence rate and provide
a limit theorem for moments in a general framework.
We also provide further insights and generalize the main result of
\cite{ARA} on robust variants of the FDH estimator, and
we provide strongly consistent and asymptotically normal estimators
$\hat\rho_x$ and $\tilde\rho_x$ of
the unknown conditional tail index $\rho_x$ involved in the limit law of
$\hat{\varphi}_1(x)$.
Moreover, when the joint density of $(X,Y)$ decreases to zero or
increases toward
infinity at a speed of
power $\beta_x>-1$ of the distance from the boundary, as is often
assumed in the literature, we answer the question of how $\rho_x$ is
linked to the data dimension $p+1$ and to the shape parameter
$\beta_x$. The quantity $\beta_x\not=0$ describes the rate at which
the density
tends to infinity (in the case $\beta_x<0$) or to $0$ (in the case
$\beta_x>0$) at the
boundary. When $\beta_x=0$, the joint density is strictly positive on the
frontier. We establish that $\rho_x=\beta_x+(p+1)$.
As an immediate
consequence, we extend the previous results of \cite{HWA,PaSW} to the
general setting where $p\geq1$
and $\beta=\beta_x$ may depend on $x$.

We propose new extreme-value-based frontier estimators $\hat\varphi^*_1(x)$,
$\tilde\varphi^*_1(x)$ and $\hat{\varphi}(x)$, which
are asymptotically normally
distributed and provide useful asymptotic confidence bands for the monotone
frontier function $\varphi(x)$.
These estimators have the advantage of not being
limited to a bi-dimensional support and
benefit from their explicit and easy formulations, which is not the
case for
estimators defined by optimization problems, such as local polynomial
estimators (see, for example, \cite{HAL}).
Their asymptotic normality is derived under quite natural and general
extreme value conditions, without Lipschitz conditions on the boundary
and without
recourse to assumptions either on the marginal distribution of $X$ or
on the
conditional distribution of $Y$ given $X=x$, as is often the case in both
statistical and econometrics literature on frontier estimation.
The study of the asymptotic properties of the different estimators considered
in the present
paper is easily carried
out by relating them to a simple dimensionless random sample and then
applying standard extreme value theory (for example, \cite{DDH,DEKK}).

Two closely related works in boundary estimation via extreme value theory
are \cite{HALL}, in which the estimation of the frontier function at a
point $x$ is based on
an increasing number of higher order statistics generated by the $Y_i$
observations
falling into a strip around $x$, and
\cite{Gij}, in which estimators are instead based on a fixed number of higher
order statistics. The main difference with the present approach is that
Hall \textit{et al.}~\cite{HALL} only focus on estimation of the
support curve of a bivariate density (that is, $p=1$)
in the case $\beta_x>1$ (that is, the decrease in density is no more
than algebraically fast),
where it is known that estimators based on an increasing number of
higher order
statistics give optimal convergence rates.
In contrast, Gijbels and Peng \cite{Gij} consider the maximum of all
$Y_i$ observations
falling into a strip around $x$
and an end-point type of estimator based on three large order
statistics of the
$Y_i$'s in the strip. This methodology is closely related and
comparable to
our estimation method using the Pickands-type estimator, but, like the
procedure of
\cite{HALL}, it is only
valid in the simple case $p=1$ and involves,
in addition to the sequence $k_n$, an extra smoothing parameter
(bandwidth of the strip) which also needs to be selected. Moreover, the
asymptotic results in \cite{Gij} are provided for densities of $(X,Y)$
decreasing as a power of the distance from the boundary, whereas the
setup in
our approach is a general one. Also, note that our transformed dimensionless
data set $(Z^{x}_1,\ldots,Z^{x}_n)$ is constructed in such a way as to
take into account the monotonicity of the
frontier (the end-point of the common distribution of the $Z^{x}_i$'s coincides
with the frontier function $\varphi(x)$), the univariate random
variables $Z^{x}_i$ do not depend on the sample size and they allow the
available results from standard extreme value theory to be easily
employed, which is
not the case for either of \cite{Gij,HALL}.

It should be clear that the monotonicity constraint on the frontier is the
main difference with most of the existing approaches in the statistical
literature. Indeed, the joint support of a random vector $(X,Y)$ is often
described in the literature as the set $\{(x,y)\mid y\leq\phi(x)\}$,
where the graph of
$\phi$ is interpreted as its upper boundary.
As a matter of fact, the function of interest, $\varphi$, in our
approach is the smallest
monotone non-decreasing function which is greater than or equal to the frontier
function $\phi$. To our knowledge, only the estimators FDH and DEA
estimate the quantity~$\varphi$. Of course,
$\phi$ coincides with $\varphi$ when the boundary curve is monotone,
but the
construction of estimators of the end-point $\phi(x)$ of the
conditional distribution of $Y$
given $X=x$ requires a smoothing procedure, which is not the case
when the distribution of $Y$ is conditioned by $X\leq x$.

We illustrate how the large-sample theory applies in practice by
carrying out some Monte Carlo experiments. Good
estimates of $\varphi(x)$ and $\rho_x$ may require
a large sample of the order of several thousand.
Theoretically selecting the
optimal extreme conditional quantiles $\hat{\varphi}_{\alpha
(k_n(x))}$ for
estimating $\varphi(x)$ and/or $\rho_x$ is a difficult question that
is worthy of future
research. Here, we suggest a simple automatic data-driven method that
provides a
reasonable choice of the sequence $\{k_n(x)\}$ for large samples.

The empirical study reveals that the simultaneous estimation of the
tail index
and of the frontier function requires large sample sizes to provide sensible
results.
The moment estimators of $\rho_x$ and of $\varphi(x)$ sometimes
provide better
estimations than the Pickands estimates and sometimes not.
When considering bias and MSE, $\hat\varphi(x)$ and $\tilde\rho_x$
provide more accurate
estimations, but when the sample size is large enough, $\hat\varphi^*_1(x)$
and $\hat\rho_x$ significantly improve and even seem to outperform
the moment
estimators. As far as the inference on $\rho_x$ is concerned, $\tilde
\rho_x$
also provides quite reliable confidence intervals, but $\hat\rho_x$ provides
more satisfactory results for sufficiently large samples.
However, when inference about the frontier function itself is
concerned, the
moment estimator provides very poor results compared with the Pickands
estimator.

On the other hand, the performance of the estimator $\tilde\varphi^*_1(x)$,
computed when $\rho_x$ is known, is quite remarkable, even compared
with the
popular FDH. The confidence intervals for $\varphi(x)$ are very
easy to
compute and have quite good coverages. In addition, the results are quite
stable with respect to the choice of the `smoothing' parameter
$k_n(x)$. As
shown in our illustrations, the estimates also have the advantage of
being robust
to extreme values. This suggests, even if $\rho_x$ is unknown, the use
of a
plug-in version of $\tilde\varphi^*_1(x)$ for making inference on
$\varphi(x)$: here, in a first step, we estimate $\rho_x$ (using the moment
estimator, unless $N_x$ is large enough), then we use the asymptotic
results for $\tilde\varphi^*_1(x)$, as if $\rho_x$ was known.
A sensible practice is not to restrict the first step to
one procedure, but rather to check that both Pickands and moment estimators
point toward similar conclusions.

\begin{appendix}\label{append}
\section*{Appendix: Proofs}

\begin{pf*}{Proof of Theorem \ref{thm1}}
Let $Z^x = Y \mathbh{1}(X\leq x)$ and $F_x (\cdot) =\{1 - F_X (x)
[1-F(\cdot|x)]\}\mathbh{1}(\cdot\geq0)$.
It can be easily seen that $\mathbb{P}(Z^x \leq y) = F_x (y)$ for any
$y\in\mathbb{R}
$. Therefore, $\{Z^x_i = Y_i \mathbh{1}(X_i \leq x)$, $i=1,\ldots ,n\}$ is an
i.i.d.~sequence of random variables with common distribution function~$F_x$. Moreover, it is easy to see that the right end-point of $F_x$
coincides with $\varphi(x)$ and that $\max_{i=1,\ldots ,n} Z^x_i$
coincides with $\hat{\varphi}_1 (x)$. Thus, assertion (i) follows
from the
Fisher--Tippett theorem. It is well known that the normalized maxima
$b^{-1}_n (\hat{\varphi}_1 (x) - \varphi(x))
\stackrel{d}{\longrightarrow} G$ (that is, $F_x$ belongs to the domain of
attraction of $G=\Psi_{\rho_x}$) if and only if
%
\begin{equation}
\label{(A.1)}
\bar{F}_x\bigl(\varphi(x) - 1/t\bigr)\in \mathit{RV} _{-\rho_x} ,
\end{equation}
where $\bar{F}_x = 1-F_x$. This
necessary and sufficient condition is equivalent to \eqref{(2.2)}. In this
case, the norming constant $b_n$ can be taken to be equal to $\varphi
(x) - \inf\{y \geq
0\mid F_x (y) \geq1-\frac{1}{n}\}=\varphi(x) -
\inf\{y \geq0\mid F(y|x) \geq1 - \frac{1}{nF_X (x)}\}$, which gives assertion
(ii). For assertion (iii), since \eqref{(A.1)} holds
and $\mathbb{E}[|Z^x|^k]=F_X(x)\mathbb{E}(Y^k|X\leq x)\leq{\varphi
(x)}^k$, it is
immediate (see \cite{RES}, Proposition~2.1) that $\lim_{n\rightarrow\infty}\mathbb{E}\{b^{-1}_n(\hat
{\varphi}_1
(x) -
\varphi(x))\}^k=(-1)^k\Gamma(1+k/\rho_x)$. Likewise, the last result follows
from~\cite{RES}, Corollary 2.3.
\end{pf*}

\begin{pf*}{Proof of Corollary \ref{(corcor)}}
Following the proof of Theorem \ref{thm1}, we can set $b_n=\varphi(x)
- F^{-1}_x(1-\frac{1}{n})$,
where $F^{-1}_x(t)=\inf\{y\in\,]0,\varphi(x)]\dvtx  F_x (y) \geq t\}$ for all
$t\in\,]0,1]$. It follows from \eqref{fxy} that $F^{-1}_x (t) =
\varphi(x) -((1-t)/\ell_x)^{1/\rho_x}$ as $t\uparrow1$
and so
$b_n=(1/n\ell_x)^{1/\rho_x}$ for all $n$ sufficiently large.
\end{pf*}

\begin{pf*}{Proof of Corollary \ref{(corscors)}}
Under the given conditions, it can be easily seen from \eqref{fxy} that
\begin{eqnarray*}
f(x,y) &=&\bigl (\varphi(x) - y\bigr)^{\rho_x-(p +
1)}\\
&&{}\times \biggl[\ell_x\rho_x(\rho_x-1)\cdots
(\rho_x-p)\frac{\partial}{\partial x^1}\varphi(x)\cdots\frac
{\partial}{\partial x^p}\varphi(x) +\mathrm{o}(1)\biggr]\qquad  \mbox{as }  y\uparrow\varphi(x) ,
\end{eqnarray*}
where the term $\mathrm{o}(1)$ depends on the partial
derivatives of $x\mapsto\ell_x$, $x\mapsto\rho_x$ and $x\mapsto
\varphi(x)$.
\end{pf*}

For the next proofs, we need the following lemma whose proof is quite
easy and
is thus omitted.

\begin{Lemma}\label{delt}
Let $Z^x_{(1)}\leq\cdots\leq Z^x_{(n)}$ be the order statistics generated
by the random variables $Z^x_1,\ldots ,Z^x_n$:
\begin{enumerate}[(iii)]
\item[(i)] If $\hat{F}_{X}(x)>0$, then
$\hat\varphi_{1-{k}/{(n\hat F_X(x))}}(x)=Z^x_{(n-k)}$
for each $k\in\{0,1,\ldots,n\hat{F}_{X}(x)-1\}$.
\item[(ii)] For any fixed integer $k\geq0$, we have
$\hat\varphi_{1-{k}/{(n\hat F_X(x))}}(x)=Z^x_{(n-k)}$ as $n\to
\infty$,
with probability $1$.
\item[(iii)] For any sequence of integers $k_n\geq0$ such
that $k_n/n\to0$ as $n\to\infty$, we have
\[
\hat\varphi_{1-{k_n}/{(n\hat F_X(x))}}(x)=Z^x_{(n-k_n)}
\qquad \mbox{as }  n\to\infty,
 \mbox{ with probability }  1.
\]
\end{enumerate}
\end{Lemma}

\begin{pf*}{Proof of Theorem \ref{Thm2.2}}
(i)
Since $\varphi(x)=F^{-1}_x(1)$ and $\hat{\varphi}_1(x)=Z^x_{(n)}$
for all $n\geq1$, we have~$(\hat{\varphi}_1(x) - \varphi(x))=(Z^x_{(n)} - F^{-1}_x(1))$. Hence,
if $b^{-1}_n (\hat{\varphi}_1 (x) -
\varphi(x))\stackrel{{d}}{\longrightarrow} G_x $, then $b^{-1}_n
(Z^x_{(n)} -
F^{-1}_x(1))$ converges to the same distribution
$G_x$.
Therefore, following \cite{vdV}, Theorem 21.18,
$b^{-1}_n (Z^x_{(n-k)}-F^{-1}_x(1))\stackrel{d}{\rightarrow} H_x$ for
any integer $k\geq0$,
where $H_x(y) = G_x(y) \sum^k_{i=0} (-\log
G(y))^i/i!$.
Finally, since $Z^x_{(n-k)}\stackrel{\mathrm{a.s.}}{=}\hat\varphi
_{1-{k}/{(n\hat F_X(x))}}(x)$
as $n\to\infty$, in view of Lemma
\ref{delt}(ii), we obtain $b^{-1}_n (\hat\varphi_{1-{k}/{(n\hat
F_X(x))}}(x)-F^{-1}_x(1))\stackrel{d}{\longrightarrow} H_x$.
{\smallskipamount=0pt
\begin{longlist}
\item[(ii)] Writing $b^{-1}_n (\hat{\varphi}_{\alpha} (x) -\varphi(x))=b^{-1}_n
(\hat{\varphi}_{\alpha} (x) -\hat{\varphi}_1 (x))+b^{-1}_n (\hat
{\varphi}_1
(x) - \varphi(x))$,
it suffices to find an appropriate sequence $\alpha=\alpha
_n\rightarrow
1$ such that $b^{-1}_n(\hat{\varphi}_{\alpha_n} (x) -\hat{\varphi
}_1 (x))\stackrel{{d}}{\longrightarrow} 0$.
Aragon \textit{et al.}~\cite{ARA} (see equation (20)) showed that
$|\hat{\varphi}_{\alpha} (x) -\hat{\varphi}_1 (x)|\leq(1-\alpha
)n\hat{F}_{X}
(x)F^{-1}_Y(1),$ with probability $1$,
for any $\alpha>0$. It thus suffices to
choose $\alpha=\alpha_n\rightarrow1$ such that $nb^{-1}_n(1-\alpha
_n)\rightarrow0$.\qed
\end{longlist}}
\noqed\end{pf*}

\begin{pf*}{Proof of Theorem \ref{thm2}}
(i) Let
$\gamma_x=-1/\rho_x$ in \eqref{(A.1)}. The Pickands \cite{Pick}
estimate of the exponent of variation $\gamma_x <0$ is then given by
$
\hat\gamma_x := (\log2)^{-1}\log\{(Z^x_{(n-k+1)} - Z^x_{(n-2k+1)})/
(Z^x_{(n-2k+1)} - Z^x_{(n-4k+1)})\}.
$
Under \eqref{(2.2)}, Condition \eqref{(A.1)} holds and so there
exists $b_n >0$ such that
$
\lim_{n\to\infty} \mathbb{P}[ b_n^{-1}(Z^x_{(n)}-\varphi(x))\le y]
= \Psi
_{-1/\gamma_x}(y).
$
Since this limit is unique only up to affine transformations, we have
\begin{eqnarray*}
\lim_{n\to\infty} \mathbb{P}\bigl[ c_n^{-1}\bigl(Z^x_{(n)}-d_n\bigr)\le y
\bigr] =
\Psi_{-1/\gamma_x}(-\gamma_x y-1)= \exp\{-(1+\gamma_x
y)^{-1/\gamma_x} \}
\end{eqnarray*}
for all $y\le0$, where $c_n= -\gamma_x b_n$ and $d_n= \varphi(x) -
b_n$. Thus, condition (1.1) from Dekkers and de Haan \cite{DDH} holds.
Therefore, $\hat
\gamma_x \stackrel{p}{\to} \gamma_x$ if $k_n\to\infty$ and
$\frac{k_n}{n}\to0$, in view of \cite{DDH}, Theorem 2.1. This gives
the weak consistency of $\hat\rho_x$ since $\hat
\gamma_x\stackrel{\mathrm{a.s.}}{=}-1/\hat\rho_x$ as $n\to\infty
$, in view of Lemma~\ref{delt}(iii).
{\smallskipamount=0pt
\begin{longlist}
\item[(ii)] Likewise, if $\frac{k_n}{n}\to0$ and $\frac{k_n}{\log\log
n}\to\infty$,
then $\hat\gamma_x \stackrel{\mathrm{a.s.}}{\longrightarrow}
\gamma_x$ via \cite{DDH},
Theorem 2.2, and so $\hat\rho_x\stackrel{\mathrm{a.s.}}{\longrightarrow}\rho_x$.

\item[(iii)] We have $U(t)= \inf\{y \ge0 \mid  \frac{1}{1-F_x(y)} \ge t\}$, which
corresponds to the inverse function $(1/(1-F_x))^{-1}(t)$. Since
$\pm t^{1-\gamma_x}U'(t) \in\Pi(A)$ with $\gamma_x=-1/\rho_x <0$,
it follows
from \cite{DDH} (see Theorem 2.3) that $\sqrt{k_n}(\hat\gamma_x -
\gamma_x)\stackrel{{d}}{\longrightarrow} {\mathcal{N}}(0,\sigma
^2(\gamma_x))$ with\vspace*{1pt}
$\sigma^2(\gamma_x)= \gamma_x^2(2^{2\gamma_x +1} +1)/\{
2(2^{\gamma_x}
-1)\log2 \}^2$ for $k_n \to\infty$ satisfying $k_n =
\mathrm{o}(n/g^{-1}(n))$, where $g(t):= t^{3- 2\gamma_x}\{
U'(t)/A(t)\}^2$.
By using the fact that
$\sqrt{k_n} (\hat\rho_x - \rho_x)\stackrel{\mathrm{a.s.}}{=}\sqrt
{k_n} (-\frac
{1}{\hat\gamma_x}
+\frac{1}{\gamma_x})$ as $n\to\infty$, in view of Lemma \ref{delt}(iii)
and applying the delta method, we conclude that
$\sqrt{k_n} (\hat\rho_x - \rho_x) \stackrel{{d}}{\longrightarrow}
{\mathcal{N}
}(0,\sigma^2(\rho_x))$
with asymptotic variance $\sigma^2(\rho_x)=\sigma^2(\gamma
_x)/\gamma_x^4$.

\item[(iv)] Under the regularity condition, we have $\pm\{
t^{-1-{1}/{\gamma_x}}F_x'(\varphi(x) - \frac{1}{t}) - \delta
F_X(x)\}\in\mathit{RV}_{-\kappa}$. The conclusion then
follows immediately from Theorem 2.5 of \cite{DDH} in
conjunction with Lemma~\ref{delt}(iii).\qed
\end{longlist}}
\noqed\end{pf*}

\begin{pf*}{Proof of Theorem \ref{thm?}}
We have, by Lemma \ref{delt}(iii), that for each $j=1,2$,
%
\begin{equation}\label{tion}
M^{(j)}_n=(1/k)\sum_{i=0}^{k-1}\bigl(\log Z^x_{(n-i)}- \log
Z^x_{(n-k)}\bigr)^j\qquad  \mbox{as }  n\to\infty ,
\mbox{ with probability } 1 ;
\end{equation}
$-1/\tilde\rho_x$ then coincides almost surely, for all $n$
large enough, with the
well-known moment estimator $\tilde\gamma_x$ (given by \cite{DEKK},
equation (1.7)) of the index defined in \eqref{(A.1)} by $\gamma
_x=-1/\rho_x$. Hence, Theorem~\ref{thm?}(i) and (ii) follow from the
weak and strong consistency of
$\tilde\gamma_x$ proved in \cite{DEKK}, Theorem 2.1. Likewise,
Theorem~\ref{thm?}(iii) follows by applying \cite{DEKK}, Corollary 3.2, in
conjunction with the delta method.
\end{pf*}

\begin{pf*}{Proof of Theorem \ref{thm3}}
(i) Under the regularity
condition, the distribution function $F_x$ of $Z^x$
has a positive derivative $F_x'(y)=F_X(x)F'(y| x)$ for all $y>0$
such that
$F'_x(\varphi(x)-\frac{1}{t})\in \mathit{RV} _{1 + {1}/{\gamma_x}}$.
Therefore, according to \cite{DDH} (see Theorem 3.1),
\[
 \sqrt{2k_n} \frac
{Z^x_{(n-k_n+1)}-F_x^{-1}(1-p_n)}{Z^x_{(n-k_n+1)}-Z^x_{(n-2k_n+1)}}
\]
 is
asymptotically normal with mean zero and variance $2^{2\gamma_x +
1}\gamma_x^2/(2^{\gamma_x}-1)^2$.
We conclude by using the facts that $F_x^{-1}(1-p_n)=\varphi_{1-{p_n}/{F_X(x)}}(x)$ and
\begin{eqnarray*}
 &&\sqrt{2k_n} \frac
{Z^x_{(n-k_n+1)}-F_x^{-1}(1-p_n)}{Z^x_{(n-k_n+1)}-Z^x_{(n-2k_n+1)}}\\
&&\quad \stackrel
{\mathrm{a.s.}}{=}\sqrt{2k_n} \frac{\hat
\varphi_{1-{(k_n-1)}/{(n\hat F_X(x))}}(x)- F_x^{-1}(1-p_n)}{\hat
\varphi_{1-{(k_n-1)}/{(n\hat F_X(x))}}(x)-
\hat\varphi_{1-{(2k_n-1)}/{(n\hat F_X(x))}}(x)}\qquad  \mbox{as }  n\to\infty .
\end{eqnarray*}
{\smallskipamount=0pt
\begin{longlist}
\item[(ii)] We have $\hat\varphi^*_1(x)\stackrel{\mathrm{a.s.}}{=}\frac
{Z^x_{(n-k_n+1)} -
Z^x_{(n-2k_n+1)}}
{2^{-\hat\gamma_x} - 1} + Z^x_{(n-k_n+1)}$ as $n\to\infty$. Following
\cite{DDH}, Theorem~3.2,
\[
\frac{\sqrt{2k_n}(\hat\varphi^*_1(x)-\varphi
(x))}{Z^x_{(n-k_n+1)} - Z^x_{(n-2k_n+1)}}
\]
 is then asymptotically
normal with mean zero and variance $3\gamma_x^2 2^{2\gamma
_x-1}/(2^{\gamma_x}-1)^6$.

\item[(iii)] Let $E_{(1)}\le\cdots\le E_{(n)}$ be the order statistics of
i.i.d.~exponential variables $E_1,\ldots, E_n$. Then, $\{
Z^x_{(n-k+1)}\}_{k=1}^n \stackrel{{d}}{=} \{U(\mathrm{e}^{E_{(n-k+1)}})\}_{k=1}^n$.
Writing $V(t):=U(\mathrm{e}^t)$, we obtain
\begin{eqnarray*}
&&\sqrt{2k_n}\biggl\{\frac{1}{2^{-\gamma_x}-1} +
\frac{Z^x_{(n-k_n+1)} -\varphi(x)}
{Z^x_{(n-k_n+1)}-Z^x_{(n-2k_n+1)}}
\biggr\}\\
&&\quad \stackrel{{d}}{=} \sqrt{2k_n}\biggl\{\frac{1}{2^{-\gamma_x}-1} +
\frac{V(E_{(n-k_n+1)}) -\varphi(x)}
{V(E_{(n-k_n+1)})-V(E_{(n-2k_n+1)})}
\biggr\}\\
&&\quad  = \biggl[-\sqrt{2k_n}\biggl\{ \frac{V(\infty)- V(\log{n}/{(2k_n)})}{V'(\log{n}/{(2k_n)})} + \frac{1}{\gamma_x}\biggr\}
 \\
&&\qquad \hspace*{4pt}{} +\sqrt{2k_n}\biggl\{ \frac{V(E_{(n-k_n+1)})-V(E_{(n-2k_n+1)})}
{2^{\gamma_x}V'(E_{(n-2k_n+1)})} -\frac{1-2^{-\gamma_x}}{\gamma
_x}\biggr\} \frac{2^{\gamma_x}}{1-2^{\gamma_x}}
\frac{V'(E_{(n-2k_n+1)})}{V'(\log{n}/{(2k_n)})}
\\
&&\hspace*{4pt}\qquad {}
- \frac{\sqrt{2k_n}}{\gamma_x} \biggl\{
\frac{V'(E_{(n-2k_n+1)})}{V'(\log{n}/{(2k_n)})} -1 -
\gamma_x\frac{V(E_{(n-k_n+1)})-V(\log{n}/{(2k_n)})}
{V'(\log{n}/{(2k_n)})}\biggr\}\biggr]\\
&&\hspace*{33pt}{}\times
\frac{V'(\log{n}/{(2k_n)})}{V(E_{(n-k_n+1)})-V(E_{(n-2k_n+1)})} .
\end{eqnarray*}
The first term on the right-hand side tends to zero as established by Dekkers
and de Haan (\cite{DDH}, Proof of Theorem 3.2). The second term
converges in distribution to ${\mathcal{N}}(0,1)\times
\frac{2^{\gamma_x}}{1-2^{\gamma_x}}$, in view of Lemma 3.1 and
\cite{DDH}, Corollary 3.1. The third term converges in probability to
$\frac{\gamma_x}{2^{\gamma_x}-1}$ by the same Corollary 3.1. This
ends the
proof of (iii), in conjunction with the fact that
\begin{eqnarray*}
&&\sqrt{2k_n} \frac{\tilde\varphi^*_1(x) - \varphi(x)}
{\hat\varphi_{1-{(k_n-1)}/{(n\hat F_X(x))}}(x)- \hat\varphi
_{1-{(2k_n-1)}/{(n\hat F_X(x))}}(x)}\\
&&\quad =\sqrt{2k_n} \biggl\{\frac{1}{2^{-\gamma_x}-1} +
\frac{Z^x_{(n-k_n+1)} -\varphi(x)}
{Z^x_{(n-k_n+1)}-Z^x_{(n-2k_n+1)}}
\biggr\}\qquad  \mbox{as }  n\to\infty ,
\end{eqnarray*}
with probability $1$.\qed
\end{longlist}}
\noqed\end{pf*}

\begin{pf*}{Proof of Theorem \ref{thm3bis}}
Write $\bar{F}_x (y) := F_X (x) [1-F(y|x)]$ and $F_x (y)
:=1 -\bar{F}_x (y)$ for all $y \geq0$.
Let $R_x (y) := - \log\{\bar{F}_x (y)\}$
for all $y\in[0, \varphi(x)[$ and let $E_{(n-k_n +1)}$ be the
statistic of order
$n-k_n +1$ generated by $n$ independent
standard exponential random variables. $Z^x_{(n-k_n +1)}$ then has
the same distribution as $R^{-1}_x [E_{(n-k_n +1)}]$, where
$
R^{-1}_x (t): = \inf\{y \geq0\mid R_x (y) \geq t\} =
\inf\{y\geq0\mid F_x(y) \geq1-\mathrm{e}^{-t}\}:= F_x^{-1} (1-\mathrm{e}^{-t}).
$
Hence,
\begin{eqnarray*}
&&Z^x_{(n-k_n +1)} - F_x^{-1} \biggl(1-\frac{k_n}{n}\biggr)\\
&&\quad  \stackrel{d}{=}
R^{-1}_x\bigl[E_{(n-k_n +1)}\bigr]-R^{-1}_x\biggl[\log
\biggl(\frac{n}{k_n}\biggr)\biggr]\\
&&\quad =\biggl[E_{(n-k_n +1)} - \log\biggl(\frac{n}{k_n}\biggr)\biggr]
(R^{-1}_x)' \biggl[\log\biggl(\frac{n}{k_n}\biggr)\biggr]\\
&&\qquad {}+ \frac
{1}{2} \biggl[E_{(n-k_n +1)} - \log
\biggl(\frac{n}{k_n}\biggr)\biggr]^2 (R^{-1}_x)'' [\delta_n] ,
\end{eqnarray*}
provided that $E_{(n-k_n +1)} \wedge\log(n/k_n)
< \delta_n < E_{(n-k_n +1)} \vee\log(n/k_n)$.
By the regularity condition~\eqref{fxy}, we have that
$ R^{-1}_x (t) =
\varphi(x) - (\mathrm{e}^{-t}/\ell_x)^{1/\gamma_x}$
for all $t$ large enough. Therefore, for
all $n$ sufficiently large,
\begin{eqnarray*}
&&\{\rho_xk^{1/2}_n/ (k_n
/n\ell_x)^{1/\rho_x} \}\bigl[Z^x_{(n-k_n +1)} -
F^{-1}_x(1-k_n/n)\bigr]\\
&&\quad \stackrel{d}{=} k^{1/2}_n\bigl[E_{(n -k_n +1)}
- \log(n/k_n)\bigr]\\
&&\qquad {} -\{k^{1/2}_n/2\rho_x\}\bigl[E_{(n-k_n+1)} -\log(n/k_n)\bigr]^2
\exp\{-[\delta_n - \log(n/k_n)]/\rho_x\} .
\end{eqnarray*}
Since $k^{1/2}_n [E_{(n-k_n +1)}-\log
(n/k_n)] \stackrel{{d}}{\rightarrow} \mathcal{N}(0,1)$
and
$|\delta_n - \log(n/k_n)|\leq|E_{(n-k_n +1)} -\break
\log(n/k_n)| \stackrel{p}{\rightarrow} 0$ as $n
\rightarrow\infty$, we obtain
$
\{\rho_xk^{1/2}_n/(k_n/n\ell_x)^{1/\rho_x}\}[Z^x_{(n-k_n +1)} - F^{-1}_x
(1-k_n/n)]\stackrel{{d}}{\longrightarrow}
\mathcal{N}(0,1)$ as $n\rightarrow\infty$.
Since $F^{-1}_x (t) = \varphi(x) - ((1-t)/\ell_x
)^{1/\rho_x}$ for all $t<1$ large enough, we have
$\varphi(x) - F^{-1}_x (1-\frac{k_n}{n}) =
(k_n/n\ell_x)^{1/\rho_x}$ for all $n$
sufficiently large. Thus,
$
\{\rho_xk^{1/2}_n/(k_n/n\ell_x)^{1/\rho_x}\}\times  [Z^x_{(n-k_n +1)}
+(k_n/n\ell_x)^{1/\rho_x} -
\varphi(x)]\stackrel{{d}}{\to} \mathcal{N}(0,1)$ as $n\to\infty$.
We conclude by using the fact that $Z^x_{(n-k_n+1)}\stackrel{\mathrm{a.s.}}{=}
\hat{\varphi}_{1-{(k_n-1)}/{(n\hat{F}_{X} (x))}}(x)$ as $n\to
\infty$.
\end{pf*}

\begin{pf*}{Proof of Theorem \ref{thm2.7}}
(i) As shown
in the proof
of Theorem \ref{thm3}(i), we have
$F'_x(\varphi(x)-\frac{1}{t})\in \mathit{RV} _{1 + 1/\gamma_x}$.
Then, by
applying Dekkers \textit{et al.}~\cite{DEKK}, Theorem 5.1, in conjunction
with \eqref{tion}, we get
\[
\sqrt{k_n}\bigl\{Z^x_{(n-k_n)}-F^{-1}_x(1-p_n)\bigr\}/M^{(1)}_nZ^x_{(n-k_n)}
\stackrel{{d}}{\longrightarrow} {\mathcal{N}}\bigl(0,V_4(-1/\gamma_x)\bigr) .
\]
The proof is completed by simply using the fact that
$F^{-1}_x(1-p_n)=\varphi_{1-{p_n}/{(F_{X} (x))}}(x)$ and
$Z^x_{(n-k_n)}\stackrel{\mathrm{a.s.}}{=}\hat{\varphi}_{1-
{k_n}/{(n\hat{F}_{X}
(x))}}(x)$ as
$n\to\infty$.
{\smallskipamount=0pt
\begin{longlist}
\item[(ii)] Since $Z^x_{(n-k_n)}\stackrel{\mathrm{a.s.}}{=}
\hat{\varphi}_{1-{k_n}/{(n\hat{F}_{X} (x))}}(x)$ and
$\tilde\gamma_x\stackrel{\mathrm{a.s.}}{=}-1/\tilde\rho_x$ as
$n\to\infty$, we have
$\hat{\varphi}(x)\stackrel{\mathrm
{a.s.}}{=}Z^x_{(n-k_n)}M^{(1)}_n(1-1/\tilde\gamma
_x)+Z^x_{(n-k_n)}$ as
$n\to\infty$.
It is then easy to see from \eqref{tion} that
$\hat{\varphi}(x)$ coincides almost surely, for all $n$ large enough,
with the
end-point estimator $\hat{x}^{*}_n$ of $F^{-1}_x(1)$
introduced by \cite{DEKK}, equation (4.8).
It is also easy to check that
$U(t)=(1/(1-F_x))^{-1}(t)$ satisfies the
conditions of \cite{DEKK}, Theorem 3.1, with $\gamma_x=-1/\rho_x<0$.
According to \cite{DEKK}, Theorem 5.2, we then have
$
\sqrt{k_n}\{\hat{x}^{*}_n-F^{-1}_x(1)\}
/M^{(1)}_nZ^x_{(n-k_n)}(1-\tilde\gamma_x)
\stackrel{{d}}{\longrightarrow} {\mathcal{N}}(0,V_5(-1/\gamma_x)),
$
which gives the desired\vspace*{-1pt} convergence in distribution of Theorem
\ref{thm2.7}(ii) since $F^{-1}_x(1)=\varphi(x)$,
$\hat{x}^{*}_n\stackrel{\mathrm{a.s.}}{=}\hat{\varphi}(x)$,
$\tilde\gamma_x\stackrel{\mathrm{a.s.}}{=}-1/\tilde\rho_x$
and $Z^x_{(n-k_n)}\stackrel{\mathrm{a.s.}}{=}
\hat{\varphi}_{1-{k_n}/{(n\hat{F}_{X} (x))}}(x)$ as $n\to\infty$.\qed
\end{longlist}}
\noqed\end{pf*}
\end{appendix}

\printhistory


\begin{thebibliography}{99}
\bibitem{ARA} Aragon, Y., Daouia, A. and Thomas-Agnan, C. (2005).
Nonparametric {f}rontier {e}stimation: {A} {c}onditional
{q}uantile-based {a}pproach. \textit{Econometric Theory} \textbf{21} 358--389.
\MR{2179542}

\bibitem{CAZ} Cazals, C., Florens, J.P. and Simar, L. (2002).
Nonparametric frontier estimation:
A robust approach. \textit{J.~Econometrics} \textbf{106} 1--25.
\MR{1875525}

\bibitem{DNN} Daouia, A. and Ruiz-Gazen, A. (2006). {Robust {n}onparametric
{f}rontier {e}stimators: {I}nfluence {f}unction and {q}ualitative
{r}obustness}. \textit{Statist. Sinica} \textbf{16} 1233--1253.
\MR{2327488}

\bibitem{DSBOOK} Daraio, C. and Simar, L. (2007).
\textit{Advanced Robust and Nonparametric Methods in Efficiency
Analysis: Methodology and Applications}. New York: Springer.

\bibitem{DDH} Dekkers, A.L.M. and de Haan, L. (1989). On the
estimation of extreme-value index and large quantiles estimation.
\textit{Ann. Statist.} \textbf{17} 1795--1832.
\MR{1026314}

\bibitem{DEKK} Dekkers, A.L.M., Einmahl, J.H.J. and de Haan, L.
(1989). A
moment estimator for the index of an extreme-value distribution.
\textit{Ann. Statist.} \textbf{17} 1833--1855.
\MR{1026315}

\bibitem{DST} Deprins, D., Simar, L. and Tulkens, H. (1984). Measuring
labor inefficiency in post offices. In \textit{The~Performance of
Public Enterprises: Concepts and Measurements} (M. Marchand, P.
Pestieau and H.~Tulkens, eds.) 243--267. Amsterdam: North-Holland.

\bibitem{Gij} Gijbels, I. and Peng, L. (2000). Estimation of a support curve
via order statistics. \textit{Extremes} \textbf{3}
251--277.
\MR{1856200}

\bibitem{HALL} Hall, P., Nussbaum, M. and Stern, S.E. (1997). On the
estimation of a support curve of indeterminate sharpness. \textit{J.
Multivariate Anal.} \textbf{62} 204--232.
\MR{1473874}

\bibitem{HAL} Hall, P., Park, B.U. and Stern, S.E. (1998). On polynomial
estimators of frontiers and boundaries. \textit{J.~Multivariate
Anal.} \textbf{66} 71--98.
\MR{1648521}

\bibitem{HAR} Hardle, W., Park, B.U. and Tsybakov, A.B. (1995). Estimation
of non-sharp support boundaries. \textit{J.~Multivariate
Anal.} \textbf{43} 205--218.
\MR{1370400}

\bibitem{HWA} Hwang, J.H., Park, B.U. and Ryu, W. (2002). Limit
theorems for
boundary function estimators. \textit{Statist. Probab. Lett.} \textbf
{59} 353--360.
\MR{1935669}

\bibitem{KST} Korostelev, A., Simar, L. and Tsybakov, A.B. (1995). Efficient
estimation of monotone boundaries. \textit{Ann. Statist.} \textbf{23} 476-489.
\MR{1332577}

\bibitem{PaSW} Park, B., Simar, L. and Weiner, C. (2000). The FDH
estimator for productivity efficiency scores:
Asymptotic properties. \textit{Econometric Theory} \textbf{16} 855--877.
\MR{1803713}

\bibitem{Pick} Pickands, J. (1975). Statistical inference using
extreme order statistics. \textit{Ann. Statist.} \textbf{3} 119--131.
\MR{0423667}

\bibitem{RES} Resnick, S.I. (1987). \textit{Extreme Values, Regular
Variation, and
Point Processes}. New York: Springer.
\MR{0900810}

\bibitem{RES2} Resnick, S.I. (2007). \textit{Heavy-Tail Phenomena:
Probabilistic and
Statistical Modeling}. New York: Springer.
\MR{2271424}

\bibitem{vdV} van der Vaart, A.W. (1998). \textit{Asymptotic Statistics.}
\textit{Cambridge Series in Statistical and Probabilistic Mathematics}
\textbf{3}. Cambridge: Cambridge Univ. Press.
\MR{1652247}

\end{thebibliography}
\end{document}